\documentclass[twoside,11pt,reqno]{amsart}

\newtheorem{propo}{Proposition}[section]

\newtheorem{lemma}[propo]{Lemma}
\newtheorem{corol}[propo]{Corollary}

\newtheorem{theo}[propo]{Theorem}

\newcommand{\bl}{\begin{lemma}\label}
\newcommand{\el}{\end{lemma}}

\newcommand{\ld}{,\ldots ,}

\newcommand{\ra}{ \rightarrow }

\newcommand{\lan}{ \langle }
\newcommand{\ran}{ \rangle }

\newcommand{\diag}{\mathop{\rm diag}\nolimits}

\newcommand{\Id}{\mathop{\rm Id}\nolimits}

\newcommand{\FF}{\mathop{\mathbb F}\nolimits}

\newcommand{\NN}{\mathbb N} 

\newcommand{\al}{\alpha}
\newcommand{\be}{\beta}

\newcommand{\ep}{\varepsilon}

\newcommand{\lam}{\lambda }

\newcommand{\up}{^{-1}}

\def\d12{{_{12}}}

\def\ag{{algebraic group }}

\def\ei{{eigenvalue }}
\def\eis{{eigenvalues }}

\def\f{{following }}

\newcommand{\med}{\medskip}

\def\ho{{homomorphism }}

\def\ii{{if and only if }}

\def\ir{{irreducible }}

\def\itf{{It follows that }}

\def\mult{{multiplicity }}

\def\Spin{{\rm Spin}}

 \newcommand{\N}{\mathbb N}

\renewcommand{\pmod}{\mod}

\begin{document}

\title[]{A sharp upper bound \\ for the size of Lusztig series}

\author{Christine Bessenrodt}

\author[A. Zalesski]{Alexandre Zalesski}

\thanks{
 E-mail:  bessen@math.uni-hannover.de (Christine Bessenrodt); alexandre.zalesski@gmail.com (Alexandre Zalesski)}

\subjclass[2000]{20C15, 20C33, 20G40}
 \keywords{Finite groups of Lie type, character theory, Lusztig series}
\maketitle

\begin{abstract} The paper is concerned
with the character theory of finite groups of Lie type. The irreducible characters of a group $G$ of Lie type are partitioned in Lusztig series. We provide a simple formula for an upper bound of the maximal size of a Lusztig series for classical groups with connected center; this is expressed for each group $G$ in terms of its Lie rank  and  defining characteristic.
When $G$ is specified as $G(q)$ and $q$ is large enough, we determine explicitly the maximum of the sizes of the Lusztig series of $G$. \end{abstract}

\section{Introduction}

Let ${\mathbf G}$ be a reductive connected algebraic group.  Let $F$ be a Frobenius (or Steinberg) endomorphism of ${\mathbf G}$ and  $G:={\mathbf G}^{F}=\{g\in G: F(g)=g\}$. Then $G$ is called a finite reductive group.

Let  ${\mathbf G}^*$ denote  the dual group  of ${\mathbf G}$, see \cite{Ca} or \cite{DM}. Then there is a Frobenius endomorphism $F^*$ of ${\mathbf G}^*$ which defines a finite group $G^*={\mathbf G}^{*F^*}=
\{g\in G^*: F^*(g)=g\}$ with $|G|=|G^*|$. The group $G^*$ is called the dual group of $G$
and plays an important role in the character theory of $G$. In particular,
the Deligne-Lusztig theory partitions the set of irreducible characters of $G$  as a disjoint union of the so called Lusztig (geometric) series ${\mathcal E}_s$, where  $s$ runs through a set of representatives
in $G^*$ of the geometric conjugacy classes of semisimple elements  of ${\mathbf G}^*$, see~\cite[13.16]{DM}. The characters from ${\mathcal E}_1$ (that is, for $s=1$) are called unipotent.

One of the questions not yet answered in the framework of the character theory of finite reductive groups is how large a Lusztig series can be. This 
has already attracted some attention in the literature, in particular, one needs to have a uniform  upper bound for $|{\mathcal E}_s|$. Liebeck and Shalev \cite[Lemma 2.1]{LS} obtained the bound $|{\mathcal E}_s|<|W|^2$, where $W$ is the Weyl group of ${\mathbf G}$, and
used this to bound the number of \ir character degrees of $G$, as well as for proving some asymptotic results. This bound has been later improved to $|{\mathcal E}_s|\leq|W|$ in \cite[Theorem 8.2]{SpZ}.

In this paper we obtain a sharp upper bound for
$\max_s|{\mathcal E}_s|$ in terms of the rank of  ${\mathbf G}$, where ${\mathbf G}$ is  a simple algebraic group of classical type with trivial center (and $s$ ranges over the semisimple elements of $G^*$). In this case $|{\mathcal E}_s|$ equals the number of unipotent characters of the group $C_{G^*}(s)$ \cite[13.23]{DM}. In fact, we compute the maximum of the number of unipotent characters of $C_{G^*}(s)$ when   ${\mathbf G}^*$ is a simple simply connected algebraic group of classical type. More precisely, we compute the maximum of $|{\mathcal E}_s|$ for $G=G(q)$ with $q$ large enough (where $q$ is the well known field parameter; usually $q>n-9$ for $q$ even and $q>n-27$ for $q$ odd, where $n$ is the rank of  ${\mathbf G}$).

To illustrate the nature of the problem, assume that $G=GL_n(q)$. Then $G\cong G^*$. If $s=1$
then $C_{G^*}(s)\cong G$.   The number of characters  in  ${\mathcal E}_1$
is well known to equal $p(n)$, the number of partitions of $n$.
One could expect that  $|{\mathcal E}_s|\leq p(n)$   for every~$s$. However, such a conjecture is false, and  the question on a sharp uniform upper bound for  $|{\mathcal E}_s|$ does not have any obvious answer. One can refine this by asking for which $s$ the number of unipotent characters of $C_{G^*}(s)$ is maximal.

In this paper we answer this question by determining the explicit value of the maximum for every classical group (for $q$ large) and describe $s$ for which the maximum is attained.  Note that it is not a priori clear at all whether the above question is feasible and can have any precise answer. The content of this paper is in computing the maximum of certain combinatorial functions to which the original problem is reduced. It is interesting and somehow surprising that the formulae we obtain for  $\max_s|{\mathcal E}_s|$ are much simpler than those available for the number of unipotent characters of $G$ and $G^*$ \cite[\S 3]{Lu}.

We expect that our results have a certain conceptual significance and will be useful for applications, in particular, for improving known upper bounds for the sum of the character degrees of $G$ (see \cite[Chapter~5]{Ko} and \cite{SpZ}).

\begin{theo}\label{un5} Let ${\mathbf G}$ be a simple \ag of rank n of adjoint type in defining characteristic $p$ and $G={\mathbf G}^F$, a finite reductive group. Then the size of a Lusztig series of $G$ does not exceed $c\cdot 5^{n/4}$ for some constant $c$ bounded as follows (and specified explicitly in the detailed results below): 

\smallskip
\begin{center}
\small{
\begin{tabular}{|c|c|c|c|c|c|
}
\hline
$A_{n}$&$C_n$&$D_n$&$B_n$&$C_n$&$D_n$\\
 \hline
&$p$ even&$p$ even&$p$ odd&$p$ odd&$p$ odd\\
\hline
$c<1.5$&$c<15$&$c<6$ &$c<95$&$c<209$&$c<44$ \\
\hline
\end{tabular}}
\smallskip
\end{center}
\end{theo}

We do not deal with the groups   ${\mathbf G}$ of exceptional Lie type in defining characteristic $p$,
and with $G={}^3D_4(q)$, as
in these cases the sizes of Lusztig series are bounded by a constant which can be easily computed.
For the other groups the constant $c$ depends on the defining characteristic $p$ of  ${\mathbf G}$, on the congruence of $n$ modulo~$4$ and, in case $D_n$, from
the choice of the Frobenius endomorphism, which defines the groups $D^+_n(q)$ or $D^-_n(q)$.
The fact that the above bound is sharp (with specified values of $c$)
can be seen from Theorem \ref{t21} below, which provides an explicit value of the maximum size of a Lusztig series for $q$ large enough
 and for each type of the group $G$. In addition, this highlights the nature of the constant $c$ and reveals that the precise
 value of $c$ in each case depends on the residue of $n$ modulo~4, a phenomenon which could not be expected in advance.

 Our starting point is a result of Bessenrodt and Ono \cite{BO}. Let $\beta(n)$ be the maximal number of the form $p(\mu):=  \Pi_j p(\mu_j)$, where $\mu=(\mu_1,\ldots, \mu_j,\ldots  )$ is a partition of $n$ and
$p(\mu_j)$ is the number of partitions of $\mu_j$.

\begin{theo}\label{bo2}  \cite{BO}
For $n=1,\ldots,7$ we have $\be(n)=p(n)=1,2,3,5,7,11,15$, respectively.
Let $\pi(n)$ denote a partition $\mu$ of $n$ such that
$\beta(n)=  p(\mu)$. Then the partition $\pi(n)$ is uniquely determined for $n\neq 7$,
whereas $\be(7)$ is attained only at the two partitions $(7)$ and $(4,3)$.

For all $n\neq 1,2,3,7$ we have the following values for $\pi(n)$ and $\beta(n)$:
\medskip

\begin{center}
{\rm Table 1}

\medskip
\begin{tabular}{|c|c|c|}\hline
$n\mod 4$& $\pi(n)$&$\beta(n)$\\
\hline
$0$&$(4^{n/4})$&$\beta(4)^{n/4}=5^{n/4}$\\
$1$&$(4^{(n-5)/4},5)$&$\beta(5)\beta(n-5)=7\cdot 5^{(n-5)/4}$\\
$2$&$(4^{(n-6)/4},6)$&$\beta(6)\beta(n-6)=11\cdot 5^{(n-6)/4}$\\
$3$&$(4^{(n-11)/4},5,6)$&$\beta(5)\beta(6)\beta(n-11)=77\cdot 5^{(n-11)/4}$\\
\hline
\end{tabular}
\smallskip
\end{center}

In particular, we always
have $\beta(n)\leq 5^{n/4}$.
\end{theo}

Using Theorem~\ref{bo2}, we obtain the following statements.

\begin{theo}\label{th1} If $G=GL_n(q)$ or $U_n(q)$ then the size of a Lusztig series does not exceed $\be(n)$,
and the bound $\be(n) = 5^{n/4}$ is attained if $4\mid n$ and $q\geq (n+1)/4$.
\end{theo}

Our results for the other classical groups are more complex.

Let $\mathbf{G}$ be a simple classical algebraic group of adjoint type and $F$ a Frobenius endomorphism such that $G=\mathbf{G}^F$ is one of the groups $SO_{2n+1}(q)$, $PCSp_{2n}(q)$, $(PSO^o)^\pm_{2n}(q)$ in the notation of \cite[Table 22.1]{MT}. Let $G^*$ be the dual group of $G$; so $G^*$ is  $Spin_{2n+1}(q)$, $Sp_{2n}(q)$ and $Spin_{2n}^\pm(q)$, respectively. For small $n$ our results, stated in Proposition \ref{p11}, are obtained by straightforward computer computations. For large $n$ a sharp upper bound for the number of unipotent characters of $C_{G^*}(s)$, and hence for the size of Lusztig series ${\mathcal E}_s$ of $G$,  when $s$ ranges over the semisimple elements of $G^*$, is provided by Theorem~\ref{t21} below.

Let $\al(m),\al^+(m),\al^-(m)$ denote the number of unipotent characters of $Sp_{2m}(q)$, $Spin^+_{2m}(q)$, $Spin^-_{2m}(q)$, respectively.

\begin{propo}\label{p11} $(1)$ Suppose that $q$ is even and $n<18$. Then $\max_s{|\mathcal E}_s|=\al(n),\al^\pm (n)$ for $G$ of type $C_n(q)$, $D^\pm_{n}(q)$, respectively.

$(2)$  Suppose that $q$ is odd and $n<32$ and $n\neq 2,4,6$ if  $G=D^-_{n}(q)$. Then $\max_s|{\mathcal E}_s|=\max_{ a}\{\al^+(a)\al(n-a)\}$, $\max_{a}\{\al(a)\al(n-a)\}$,
$\max_a\{\al^+(a)\al^+ (n-a)\}$, $\max_a\{\al^-(a)\al^+ (n-a)\}$ for $G$ of type $B_n(q)$, $C_n(q)$, $D^+_{n}(q)$, $D^-_{n}(q)$, respectively, where a ranges between $0$ and n. The explicit value of the maximum in each case is given by Tables $3,4$.
\end{propo}

\begin{theo}\label{t21}
 Let $G^*\in \{Spin_{2n+1}(q); ~Sp_{2n}(q); ~Spin^\pm_{2n}(q)\}$. For a semisimple element $s\in G^*$ let $ {\mathcal E}_s$ denote the Lusztig series of irreducible characters of~$G$.

\medskip
$(1)$  Let $G^*=Sp_{2n}(q)$, $q$ even, or $Spin_{2n}(q)$, $q$ even. 
For $n\geq 18$, we have $|{\mathcal E}_s|\leq f(n)$, where

\medskip
$f(n)=\begin{cases}
\al(16)\beta(n-16)=8988\cdot 5^{(n-16)/4}& \text{if }  ~n \equiv 0 \mod 4,\\
\al(15)\beta(n-15)=66396\cdot 5^{(n-21)/4}& \text{if }  ~n \equiv 1 \mod 4,\\
\al(14)\beta(n-14)=4020\cdot 5^{(n-14)/4}& \text{if }  ~n \equiv 2 \mod 4,\\
\al(15)\beta(n-15)=6036\cdot 5^{(n-15)/4}& \text{if }  ~n \equiv 3 \mod 4.
\end{cases}$

\medskip
$(2)$ Let $G^*=Spin^\pm_{2n}(q)$, $q$ even. 
For $n\geq 18$, we have  $|{\mathcal E}_s|\leq f^\pm(n)$, where

\medskip
$f^+(n)=\begin{cases}
\al^+ (16)\beta(n-16)=4110\cdot 5^{(n-16)/4}& \text{if }  ~n \equiv 0 \mod 4,\\
\al^+ (17)\beta(n-17)=6007\cdot 5^{(n-17)/4}& \text{if }  ~n \equiv 1 \mod 4,\\
\al^+ (14)\beta(n-14)=1836\cdot 5^{(n-14)/4}& \text{if }  ~n \equiv 2 \mod 4,\\
\al^+ (15)\beta(n-15)=2730\cdot 5^{(n-15)/4}& \text{if }  ~n \equiv 3 \mod 4,
\end{cases}$

\medskip
and

\medskip
$f^-(n)=\begin{cases}
\al^- (16)\beta(n-16)=4066\cdot 5^{(n-16)/4}& \text{if }  ~n \equiv 0 \mod 4,\\
\al^- (17)\beta(n-17)=6007\cdot 5^{(n-17)/4}& \text{if }  ~n \equiv 1 \mod 4,\\
\al^- (14)\beta(n-14)=1806\cdot 5^{(n-14)/4}& \text{if }  ~n \equiv 2 \mod 4,\\
\al^- (15)\beta(n-15)=2730\cdot 5^{(n-15)/4}& \text{if }  ~n \equiv 3 \mod 4.
\end{cases}$
\med

\medskip
$(3)$   Let  $G^*=Sp_{2n}(q)$, $q$ odd. 
For $n\geq 32$, we have $|{\mathcal E}_s| \leq \tau(n)$, where

\medskip
$\tau(n)=\begin{cases}
\al(14)\al(14)\beta(n-28)=16160400\cdot 5^{(n-28)/4}&\text{if }  ~n \equiv 0 \mod 4,\\
\al(15)\al(14)\beta(n-29)=24264720\cdot 5^{(n-29)/4}&\text{if }  ~n \equiv 1 \mod 4,\\
\al(15)\al(15)\beta(n-30)=36433296\cdot 5^{(n-30)/4}&\text{if }  ~n \equiv 2 \mod 4,\\
\al(16)\al(15)\beta(n-31)=54251568\cdot 5^{(n-31)/4}&\text{if }  ~n \equiv 3 \mod 4.
\end{cases}$

\smallskip

\medskip
$(4)$  Let  $s\in G^*=Spin_{2n+1}(q)$, $q$ odd. 
For $n\geq 32$, we have $|{\mathcal E}_s| \leq \theta(n)$, where

\medskip
$\theta(n)=\begin{cases}
\al(16)\al^+(16)\beta(n-32)=36940680\cdot 5^{(n-32)/4}&\text{if }  ~n \equiv 0 \mod 4,\\
\al(15)\al^+(14)\beta(n-29)=11082096\cdot 5^{(n-29)/4}&\text{if }  ~n \equiv 1 \mod 4,\\
\al(14)\al^+(16)\beta(n-30)=16522200\cdot 5^{(n-30)/4}&\text{if }  ~n \equiv 2 \mod 4,\\
\al(15)\al^+(16)\beta(n-31)=24807960\cdot 5^{(n-31)/4}&\text{if }  ~n \equiv 3 \mod 4.
\end{cases}$

\medskip
 $(5)$  Let 
 $G^*=Spin_{2n}^\pm (q)$,  $q$ odd. 
For $n\geq 32$, we have $|{\mathcal E}_s| \leq \theta^\pm (n)$, where

\medskip
$\theta^+(n)=\begin{cases}
\al^+(16)\al^+(16)\beta(n-32)=16892100\cdot 5^{(n-32)/4}&\text{if }  ~n \equiv 0 \mod 4,\\
\al^+(15)\al^+(14)\beta(n-29)=5012280\cdot 5^{(n-29)/4}&\text{if }  ~n \equiv 1 \mod 4,\\
\al^+(16)\al^+(14)\beta(n-30)= 7545960\cdot 5^{(n-30)/4}&\text{if }  ~n \equiv 2 \mod 4,\\
\al^+(16)\al^+(15)\beta(n-31)=11220300\cdot 5^{(n-31)/4}&\text{if }  ~n \equiv 3 \mod 4,
\end{cases}$

and

$\theta^-(n)=\begin{cases}
\al^-(16)\al^+(16)\beta(n-32)=16711260\cdot 5^{(n-32)/4}& \text{if }  ~n \equiv 0 \mod 4,\\
\al^-(15)\al^+(14)\beta(n-29)=5012280 \cdot 5^{(n-29)/4}& \text{if }  ~n \equiv 1 \mod 4,\\
\al^-(16)\al^+(14)\beta(n-30)= 7465176\cdot 5^{(n-30)/4}& \text{if }  ~n \equiv 2 \mod 4,\\
\al^-(15)\al^+(16)\beta(n-31)=  11220300 \cdot 5^{(n-31)/4}& \text{if }  ~n \equiv 3 \mod 4.
\end{cases}$

\medskip In addition, the detailed bounds given in $(1)$ - $(5)$  are attained if $q> n-9$ if $q$ is even, and $q>n-27$ if $q$ is odd.
  \end{theo}

Note that Lusztig series were  originally defined only for groups with connected center, but later
this notion has been extended to arbitrary connected reductive groups so that, again,  the size of a Lusztig series equals the number of unipotent characters of $C_{G^*}(s)$, see \cite[Theorem 13.23]{DM}.
Then the number of unipotent characters of $C_{G^*}(s)$ does not exceed  $|C_{G^*}(s): C_{G^*}(s)^0|\cdot \nu(C_{G^*}(s)^0)$,
where $C_{G^*}(s)^0=(C_{{\mathbf G}^*}(s)^0))^F$ and $\nu(C_{G^*}(s)^0)$ is the number of unipotent characters of~$C_{G^*}(s)^0$. The index $|C_{G^*}(s): C_{G^*}(s)^0|$ does not exceed $r+1$ for groups of type $A_{r}$
and 4 for the other simple groups \cite[Ch. II, Corollary 4.4]{StSp}. So we can replace $c$ by $(r+1)c$
for the $A_r$-case and $4c$ for the other groups (in fact, the latter is needed only for $q$ odd). However,
these bounds may not be sharp.

Our strategy can be outlined as follows. The simplest case is where $G=GL_n(q)$ or $U_n(q)$; here we show
(Section 3) that $|\mathcal{E}_s|\leq \be(n)$ and the bound is attained for $q$  large enough. For  the other classical groups this bound is  valid only if $\pm 1$ are not \eis of $s$ on the natural $\FF_qG^*$-module $V$ for $G^*$ (Lemma~\ref{ei1}). 
Suppose first that $G^*=Sp_{2n}(q)$, $q$ odd, 
and $k,l$ are the multiplicities of the \eis 1 and $-1$, respectively, of $s$ on $V$. Then we show that $|\mathcal{E}_s|\leq \max \al(k/2)\al(l/2)\be(m)$, where $\frac{(k+l)}{2}+m=n$. This reduces the problem to computing the above maximum, and next we show that the bound is attained for some $s$ if $q$ is large enough. If $G^*$ is an orthogonal group then we have a similar reduction with $\al(k/2)\al(l/2)$ to be replaced by $\al^\pm(k/2)\al^{\pm}(l/2)$ if $\dim V$ is even,  and $\al(k/2)\al^{\pm}(l/2)$ if $\dim V$ is odd, with a certain choice of the signs. If $q$ is even then we argue similarly. The maximum of the products in question is computed in Section~5. The proof of Theorem \ref{t21} occupies Sections 6,7, for $q$ even and odd, respectively.

\smallskip

{\it  Notation}.
 The size of a finite set $S$ is denoted by $|S|$. Also, we write $|g|$ for the order of a group element,
 which does not lead to a confusion.  For a group $G$ we denote by   $Z(G)$ the center of $G$, and  by $C_G(S)$ the centralizer of a subset $S$ of $G$ in $G$. We use this notation also in the situation where $V$ is a set on which $G$ acts by permutations or  a vector space on which $G$ acts by linear transformations. So $C_V(S)=\{v\in V: sv=v$ for all $s\in S\}$. For $S\subset G$ we denote by $\lan S\ran$ the subgroup generated by $S$.

 $\Id_n$ is the identity $(n\times n)$-matrix. By $\diag(x_1\ld x_n)$ we denote the diagonal matrix with subsequent diagonal entries $x_1\ld x_n$. A similar notation is used for a block-diagonal matrix.

We denote by $\NN$  the set of natural numbers. For  $n\in\NN$, $p(n)$ denotes the number of partitions of $n$;
for a partition  $\mu=(\mu_1,\ldots,\mu_t)$ we set $p(\mu)=\Pi_{j=1}^t p(\mu_j)$.
We then set $\beta(n)=\max p(\mu)$, where the maximum is taken over all partitions $\mu$ of $n$.

By $\FF_q$ we denote the field of $q$ elements.
If $K$ is a field then $K^\times$ denotes the multiplicative group of $K$
and  $\overline{K} $  an algebraic closure  of $K$.

All vector spaces considered in the paper are of finite dimension.
By $GL(V)$ we denote the group of all invertible linear transformations of a vector space~$V$.
If the ground field $K$ is not algebraically closed, and $s\in GL(V)$ is a semisimple element, then the natural analog of eigenspaces are homogeneous components of $s$ on $V$; these are the sum of all minimal non-zero $K\lan s\ran$-submodules of $V$ isomorphic to each other.
If $s$ has a single homogeneous component on $V$, we say that  $s$ is homogeneous.

For an algebraic group $\mathbf{G}$ we denote by $\mathbf{G}^0$ the connected component of the identity of $\mathbf{G}$.
We use $F$ to denote a Frobenius endomorphism of an algebraic group, and we often use it for different algebraic groups.
We usually write $G$ for  $\mathbf{G}^F=\{g\in \mathbf{G}: F(g)=g\}$. If $\mathbf{G}$
is  connected  reductive, we call $G=\mathbf{G}^F$ a finite reductive group.  For a finite reductive group $G$ we denote by  $\nu(G)$ the number of unipotent characters of $G$.
See Section 2.2 for more details.
As mentioned in the introduction,
$\al(m),\al^+(m),\al^-(m)$ stands for the number of unipotent characters of $Sp_{2m}(q)$, $Spin^+_{2m}(q)$, $Spin^-_{2m}(q)$, respectively.
By $\mathbf{G}^*$ and $G^*$ we denote the dual groups of a reductive \ag $\mathbf{G}$
and of a finite reductive group $G$, respectively.

Our notation for classical groups is standard, except for  the special orthogonal groups of even characteristic; following \cite{MT}, we denote by $SO^\pm_{2n}(q)$ and $SO_{2n}(\overline{\FF}_q)$
with $q$ even the subgroup of index 2 in the full orthogonal group  $O^\pm_{2n}(q)$ and $O_{2n}(\overline{\FF}_q)$, respectively. (The advantage of this is that certain results can be stated uniformly for $q$ odd and even.) In addition, dealing with   the groups $SO_{2n+1}(q)$ we assume that  $q$ is odd as  $SO_{2n+1}(q)\cong Sp_{2n}(q)$ whenever $q$ is even.

We expect a reader to be familiar with the geometry of classical groups; most necessary facts can be found in \cite[Ch. 2]{KL}. In particular, for the notion of Witt defect of an orthogonal space see \cite[p. 28]{KL}. Nonetheless we recall a few notions from this area.

An orthogonal space means a vector space $V$ of finite dimension over a field~$K$, say, endowed with a non-degenerate symmetric bilinear form $f(v_1,v_2)\in K$ for $v_1,v_2\in V$, and if the characteristic of $K$ equals 2 then the form is additionally assumed to be alternating (that is, $f(v,v)=0$ for $v\in V$) and non-defective \cite[Ch. I, \S 16]{Di}. The full orthogonal group is denoted by $O(V)$. The spinor group of an orthogonal space $V$ (we call it the full spinor group of $V$)  is defined in terms of the Clifford algebra of $V$ \cite[Ch. II, \S 7]{Di}; this yields the notion of spinor norm, which defines the subgroup $\Omega(V)$ of $O(V)$ formed by elements of spinor norm 1. In particular, if the ground field is of characteristic 2,
we have $\Omega(V)=SO(V)$ by convention. We use $Spin(V) $ to denote the preimage of $\Omega(V)$ in the spinor group of $V$ under the natural projection of it onto $O(V)$, see loc.cit.

\section{Preliminaries}

For later considerations we will need the explicit formulae for $\beta(n)$ from \cite{BO} which we now recall.

\subsection{Some properties of the function $\be(n)$.}

From Theorem~\ref{bo2}
we deduce a number of properties of the numbers $\be(n)$.

\begin{lemma}\label{bo1}
Let $0<k\leq n $ be integers. Then $\beta(k)<\beta(n)$ for $k<n$ and
 $\beta(k)\beta(n)\leq \beta(k+n).$ More precisely,
if $k\neq 1,2,3,7$ then

$$\frac{\beta(k+n)}{\beta(k)\beta(n)}=\begin{cases}
1& \text{if }  ~k\equiv 0~\text{ or }~n \equiv 0 \mod 4;\cr
\frac{55}{49} & \text{if }  ~k\equiv 1~and~n \equiv 1 \mod 4;\cr
 1& \text{if }  ~k\equiv 1~and~n \equiv 2
 ~ \text{ or }~k\equiv 2~and~n \equiv 1 \mod 4;
\cr
\frac{625}{539} & \text{if }  ~k\equiv 1~and~n \equiv 3
~\text{ or } ~k\equiv 3~and~n \equiv 1 \mod 4;
\cr
\frac{125}{121} & \text{if }  ~k\equiv 2~and~n \equiv 2 \mod 4;\cr
\frac{125}{77} & \text{if }  ~k\equiv 2~and~n \equiv 3
~\text{ or } ~k\equiv 3~and~n \equiv 2 \mod 4;\cr
\frac{625}{539} & \text{if }  ~k\equiv 3~and~n \equiv 3 \mod 4.
\end{cases}$$

\med
For $k=1,2,3,7$ and $n>7$ the values of $\frac{\be(n+k)}{\be(k)\be(n)}$ are as follows

\med
\begin{center}
\begin{tabular}{|l|c|c|c|c|c}
\hline
$n$&$k=1$&$k=2$&$k=3$&$k=7$
\cr
\hline
&&&&\\[-10pt]
$n \equiv 0 \mod 4$& $\frac{7}{5}$&$\frac{11}{10}$&$\frac{77}{75}$&$\frac{77}{75}$\\[1pt]
\hline
&&&&\\[-10pt]
$n \equiv 1 \mod 4$&$\frac{11}{5} $&$\frac{11}{10}$ &$\frac{25}{21}$& $\frac{25}{21}$ \\[1pt]
\hline
&&&&\\[-10pt]
$n \equiv 2 \mod 4$& $\frac{7}{5}$&$\frac{25}{22}$&$\frac{35}{33}$&$\frac{35}{33     }$ \\[1pt]
\hline
&&&&\\[-10pt]
$n \equiv 3 \mod 4$ &$\frac{125}{77}$&$\frac{25}{22}$&$\frac{25}{21}$&$\frac{25}{21}$ \\[1pt]
\hline\end{tabular}
\end{center}
\end{lemma}

Proof. Straightforward computations using Theorem~\ref{bo2}. $\Box$

\med

\begin{lemma}\label{bo3}
Let $n>2$ be an integer.
\begin{enumerate}

\item
If $~n$ is even, then $\beta(n/2)\leq 3\beta(n-3)$
and, for $n>4$, we have  $\beta(n/2)\leq 7\beta(n-5)$.
\item
We have $ 5^{(n-3)/4}< \be(n).$
\end{enumerate}
\end{lemma}

Proof. (1) If $n/2\leq n-5$ then $\beta(n/2)\leq \beta(n-5)<\be(n-3)$ by Lemma~\ref{bo1}. Otherwise $n<10$ and the claim follows by inspection.

(2) We have $5^{ (n-3)/4}<5^{ (n-2)/4}<5^{ (n-1)/4}<5^{n/4}$, and $5^{(n-i)/4}=\be(n-i)$ if $4|(n-i)$ with $i=0,1,2,3$. If $i>0$ then $\be(n-i)<\be(n)$ by Lemma~\ref{bo1}.
If $i=0$ then $\be(n)=5^{n/4}>5^{ (n-3)/4}$, whence the result.  $\Box$

\med
For the use in later sections we record the following lemma:

\bl{nk8} Let $n\in\NN$ be even and $n> 6$. Set $\beta'(n)=\max_{a~ odd } \,\beta (a)\beta (n-a)$. Then we have $\beta'(n)=\beta(5)\beta(n-5)$. Explicitly, we have
$$\beta'(n)=\begin{cases}
7^2\cdot 5^{(n-10)/4}  & \text{if } ~n \equiv 2 \mod 4,
\\
7^2\cdot 11\cdot 5^{(n-16)/4} & \text{if }  ~n\equiv 0 \mod 4  \text{ and } n\geq 16 ,
\end{cases}$$
and $\beta'(8)=7\cdot 3$, $\beta'(12)=7\cdot 15$.
In addition,   $\beta'(6)=9$, $\beta'(4)=3$, $\beta'(2)=1$.
\el

Proof. The additional statement follows by inspection (see Table 2). Let $i,j\in \N$ be odd with $i+j=n$, $i\le j$. Since $n\ge 10$ we have $j\ge 5$. If $j\equiv 3 \mod 4$ and $j>7$, then $\beta(j-5)=\beta(6)\beta(j-11)$ by Theorem~\ref{bo2}. Hence Theorem~\ref{bo2} and Lemma~\ref{bo1}
imply for $j\ne 7$: $\beta(i)\beta(j)=\beta(i)\beta(5)\beta(j-5)\le \beta(5)\beta(n-5)$.
For $j=7$ and $i\in \{3,7\}$, we have $\beta(3)\beta(7)= 3\cdot 15 < 7^2 =\beta(5)^2$,
and $\beta(7)^2=15^2 < 7\cdot 35 =\beta(5)\beta(9)$. So in any case $\beta'(n)=\beta(5)\beta(n-5)$.

Hence, applying Theorem~\ref{bo2} we obtain  the formulae for $\beta'(n)$ stated above. $\Box$

\subsection{Unipotent characters}

Let ${\mathbf G}$ be a connected reductive  algebraic  group with Frobenius endomorphism $F$. For a precise definition of it we refer to \cite[p. 31]{Ca} or \cite[Section 2.1]{MT}. (Some authors use the terms "Frobenius map" or "Steinberg endomorphism".) If ${\mathbf G}$ is simple then an algebraic group endomorphism  $F:{\mathbf G}\ra {\mathbf G}$ is Frobenius \ii   the subgroup ${\mathbf G}^F=\{g\in {\mathbf G}: F(g)=g\}$ is finite \cite[Theorem 21.5]{MT}.  Groups ${\mathbf G}^F$   are called {\it finite reductive groups} \cite[p. XIII]{CE} or \cite[\S 4.4]{Ca}. (The term "finite groups of Lie type" is also in use in the literature, see
\cite[p. 31]{Ca}.) Thus, a finite reductive group is determined by the pair $({\mathbf G},F)$, a connected reductive algebraic group ${\mathbf G}$ and a Frobenius endomorphism $F$ of it.

\medskip
As shortly mentioned in the introduction, for every finite reductive group $G={\mathbf G}^F$ the
 Deligne-Lusztig theory   partitions  the set of irreducible characters of $G$ as a disjoint union of the  Lusztig (geometric) series ${\mathcal E}_s$, where  $s$ runs through a set of representatives
of the  classes of semisimple elements  of $G^*$ that are conjugate
in~${\mathbf G}^*$, see~\cite[13.16]{DM}. The characters in ${\mathcal E}_1$ (that is, for $s=1$) are called {\it unipotent}. Note that the geometric series can be further refined to rational series parameterized by the conjugacy classes of semisimple elements in $G^*$; if ${\mathbf G}$ has connected center (assumed in this paper) then
the geometric and rational series coincide \cite[p. 107]{DM}.

We  emphasize that the Lusztig series (and hence unipotent characters) of a finite reductive group cannot be defined in terms of ${\mathbf G}^F$ as an abstract group.   One observes that a given finite group of Lie type can be obtained as ${\mathbf G}^F$ from different pairs ${\mathbf G},F$.   A typical example is as follows.
Given a pair $({\mathbf G},F)$, set ${\mathbf H}$ to be the direct product ${\mathbf G}\times \cdots\times{\mathbf G} $ of $m$ copies of ${\mathbf G}$, and then define a
Frobenius endomorphism $F'$ of ${\mathbf H}$ as a mapping sending an element $(g_1\ld g_m)$ with $g_1\ld g_m\in {\mathbf G}$ to $(F(g_m), g_1\ld g_{m-1})$. Then  $F'(g_1\ld g_m)=(g_1\ld g_m)$ implies $F(g_m)=g_1=g_2=\cdots =g_m$, so
${\mathbf H}^{F'}=\{(g\ld g): g\in {\mathbf G}^F\}\cong {\mathbf G}^F$. In fact,
the general case reduces to the above example, see \cite[p. 380]{Ca} where it is stated that one can assume ${\mathbf H}$ to be simple (if so is ${\mathbf G}$), that is, $m=1$.

\bl{f13} \cite[p.~112]{DM} Let $G={\mathbf G}^F$ be a finite reductive group and $s\in G^*$ a semisimple element. Suppose that $C_{{\mathbf G}^* }(s)$ is connected. Then $|{\mathcal E}_s|=\nu(C_{G^*}(s))$, the number of unipotent characters of $C_{G^*}(s)$. \el

To be rigorous, we emphasize that  $C_{G^*}(s)=C_{{\mathbf G}^*}(s)^F$ is a finite reductive group.

Lemma~\ref{f13} reduces the computation of the sizes of Lusztig series to the  computation of the number of unipotent characters, and our results in fact give sharp upper bounds for  $\nu(C_{G^*}(s))$ when $s$ ranges over semisimple elements of $G^*$. (Note that, if $C_{\mathbf G}(s)$ is not connected, one can extend the notion of a unipotent character so
that Lemma~\ref{f13} remains valid, see \cite[p. 112]{DM}. However, in full generality the problem of computing sharp upper bounds is more complex.)

\med
For what follows it is essential to decide whether $C_{\mathbf G}(s)$ is a connected reductive  group if so is ${\mathbf G}$
and $s\in {\mathbf G}$ is a semisimple element. There are the following  criteria
for connectivity:

\bl{gf4}
The group $C_{\mathbf G}(s)$ is connected for all semisimple elements $s$ of ${\mathbf G}$ if  one of the \f holds:

$(1)$ The center of ${\mathbf G}^*$ is connected.

 $(2)$  ${\mathbf G}$ is semisimple and simply connected.

 $(3)$  ${\mathbf G}=SO({\mathbf V})$, where ${\mathbf V}$ is an orthogonal space over $\overline{\FF}_q$,
 and the \mult of the \ei $1$ or the \ei $-1$ of $s$ on ${\mathbf V}$ is at most $1$.
\el

Proof. See \cite[13.15]{DM} for (1), and \cite[Ch. E, 3.9] {StSp} for (2).

(3) If $\dim V$  is even then the \mult of the \ei $1$ as well as the \ei $-1$ of $s$ on ${\mathbf V}$ is known to be even (see Lemma~\ref{ei2} below for a proof), so (3) follows from \cite[Lemma 2.2]{Z14} in this case.
Now suppose that  $\dim V$  is odd, so, by our convention, $q$ is odd.
Let ${\mathbf V}_1$ be the 1-eigenspace of $s$
on ${\mathbf V}$. Then $\dim {\mathbf V}_1=1$ (Lemma~\ref{ei2}), so $C_{\mathbf G}(s)$ is contained in the stabilizer of ${\mathbf V}_1$ in ${\mathbf G}$.
With respect to a suitable basis of ${\mathbf V}$, the latter
can be written as $\{\diag(\det g,g): g\in O({\mathbf V}^\perp_1)\}$.  Then $C_{\mathbf G}(s)$ is contained in
the group $\diag(\pm 1,C_{O({\mathbf V}^\perp_1)}(s'))$, where $s'$ is the restriction of $s$ to $   {\mathbf V}^\perp_1$. As $\dim {\mathbf V}^\perp_1$ is even, $C_{O({\mathbf V}^\perp_1)}(s')$ is connected by the above,
and hence is contained in $SO({\mathbf V}^\perp_1)$ (or see the proof of \cite[Lemma 2.1]{Z14}). Then $C_{\mathbf G}(s)=\diag(1,C_{O({\mathbf V}^\perp_1)}(s')) $, whence the result. $\Box$

\med Remark. There is an inaccuracy in the statement of \cite[Lemma 2.2]{Z14}, where "Let ${\mathbf G}=SO({\mathbf V})$" is to be replaced by "Let ${\mathbf G}=SO({\mathbf V})$ if $q$ is odd and  $\Omega({\mathbf V})$ if $q$ is even"  with no change of the proof.

\med
Thus, if Lemma~\ref{gf4} applies then $C_G(s)$ is a finite reductive group.
For the notion of a simply connected semisimple algebraic  group see for instance \cite[9.14]{MT} or \cite[p. 25]{Ca};
if ${\mathbf G}$  is of adjoint type then $G^*$ is   simply connected.
Classical algebraic groups of adjoint and simply connected type can be described in terms of their
traditional definition, see \cite[Table 9.2]{MT} or \cite[p. 40]{Ca}.

\bl{f14} Let $H=SO(V)$, where $\dim V$ is even,  and let $s\in H$ be a semisimple element.

$(1)$ Suppose that either  $1$ or $-1$ is not an \ei of $s$ on $V$. Then  $C_H(s)$ is a finite reductive group. In particular, this is the case if $q$ is even.

$(2)$ Suppose that   neither $1$ nor $-1$ is an \ei of $s$ on $V$. Then \hbox{$C_{O(V)}(s)\subset H$}.
\el

Proof. (1) Let ${\mathbf V}=V\otimes \overline{\FF}_q$ be an orthogonal space defined with the same Gram matrix as $V$.  It is well known that ${\mathbf G}=SO({\mathbf V})$ is an \ag and $SO(V)=SO({\mathbf V})^F$ for some Frobenius morphism $F: {\mathbf G}\ra {\mathbf G}$.
By \cite[Lemma~2.2(2)]{Z14},
the group $C_{{\mathbf G}}(s)$ is connected. As $C_G(s)=C_{{\mathbf G}}(s)^F$, the claim follows.

(2) See \cite[Lemma 2.1]{Z14}. $\Box$

\med
Let $\mathbf{V}$ be an orthogonal  space over $\overline{\FF}_q$.
The  group $SO(\mathbf{V})$  is a simple algebraic group, however,  $SO(\mathbf{V})$  is not simply connected. Slightly abusing notation, we denote the  simply connected covering of it by $Spin(\mathbf{V})$; this is the preimage of $SO(\mathbf{V})$ in the full spinor group of $\mathbf{V}$. So  $Spin(\mathbf{V})$ is a simply connected simple  algebraic group, and  there is a surjective algebraic group \ho
$\eta:Spin(\mathbf{V})\ra SO(\mathbf{V})$ (see \cite[p. 228]{CE}). If $q$ is even then $\eta$ 
is an isomorphism of the underlying abstract groups.

Let $h: {\mathbf G}\ra {\mathbf H}$ be a surjective  \ho of connected algebraic groups with central kernel
(that is, an isogeny), defined over $\FF_q$, and let $F$ be a Frobenius endomorphism of~${\mathbf G}$.
If ker$\,h$ is $F$-stable, one defines the action of $F$ on ${\mathbf H}$ by $F(h(g))=h(F(g))$ for $g\in {\mathbf G}$. 
Set $H:={\mathbf H}^F$. With these notations we have

\bl{f11} \cite[13.20]{DM} Let $\nu(G)$, $\nu(H)$ be the number of unipotent characters of $G,H$, respectively. Then $\nu(G)=\nu(H)$.\el

For instance, if $G=GL_n(q)$ and $H=PGL_n(q)$, or   $G=U_n(q)$ and $H=PU_n(q)$, then the lemma applies. Moreover,  $\nu(SL_n(q))=\nu(PGL_n(q))$ as $PSL_n(\overline{\FF}_q)^F=PGL_n(q)$ for $n>1$, see \cite[p. 39]{Ca}.

Lemma~\ref{f11} allows us to ignore the case where ${\mathbf G}^*=Spin_{2n+1}(\overline{\FF}_q)$ with $q$
even. Indeed, in this case there exists an isogeny $h:Spin_{2n+1}(\overline{\FF}_q)\ra {\mathbf H}:=Sp_{2n}(\overline{\FF}_q)$, which also yields an isogeny $C_{{\mathbf G}^*}(s)\ra C_{{\mathbf H}}(h(s))$.
Therefore,  by Lemma~\ref{f11}, we have $\nu(C_{G^*}(s))=\nu(C_{H}(h(s)))$, where $H={\mathbf H}^F=Sp_{2n}(q)$. So it suffices to compute the maximum of $\nu(C_{H}(s))$ over semisimple elements $s\in H$.

For the group $SO_{2n}(\overline{\FF}_q)$ there are Frobenius endomorphisms for which
$SO_{2n}(\overline{\FF}_q)^F$ coincides with $SO_{2n}^+(q)$ or $SO_{2n}^-(q)$. (If $n=4$ there is one more type of Frobenius endomorphisms which yields the "triality group" ${}^3D_4(q)$; this is not considered in this paper.) Here $SO^+_{2n}(q)$ and $SO^-_{2n}(q)$   are special orthogonal groups $SO(V)$, where $V$ is an orthogonal space of Witt defect 0 and 1, respectively, with $\dim V=2n$.

There exists a Frobenius endomorphism $F$, say, of $Spin(\mathbf{V})$ compatible with the natural mapping $\eta:Spin(\mathbf{V})\ra  SO(\mathbf{V})$ in the sense that $\eta(F(h))=F(\eta(h)) $. Then we set $Spin(V)= Spin(\mathbf{V})^F$. If $q$ is odd, then $\eta(Spin(V))=\Omega(V)\neq SO(V)$.  Nonetheless, by Lemma~\ref{f11}, we have

\bl{f12}  If $q$ is odd, then $\nu(Spin(V))=\nu(SO(V))$.\el

\begin{lemma}\label{55n} Let ${\mathbf G}=Spin({\mathbf V})$, and
let $s\in {\mathbf G}$ be a semisimple element. Let $\eta:{\mathbf G}\ra SO({\mathbf V})$
be the natural projection. Let ${\mathbf W}_1$, ${\mathbf W}_2$ be the $1$- and $-1$-eigenspaces of $\eta(s)$ on ${\mathbf V}$,
and ${\mathbf W}_3=({\mathbf W}_1+{\mathbf W}_2)^\perp$. Then $\eta(C_{{\mathbf G}}(s))=SO(\mathbf{W}_1)\times SO(\mathbf{W}_2)\times C_{SO(\mathbf{W}_3)}(s')$, where $s'\in SO({\mathbf W}_3)$ is the restriction of $\eta(s)$ to ${\mathbf W}_3$. (If ${\mathbf W}_1=0$ or ${\mathbf W}_2=0$ then the respective multiple is to be dropped.)\el

Proof. Clearly,  $\eta(C_{{\mathbf G}}(s))$ stabilizes ${\mathbf W}_1$ and ${\mathbf W}_2$, and hence also ${\mathbf W}_3$. \itf $\eta(C_{{\mathbf G}}(s))\subset O(\mathbf{W}_1)\times O(\mathbf{W}_2)\times O(\mathbf{W}_3)$. By Lemma~\ref{gf4} and the comments after Lemma~\ref{f14}, the group  $C_{{\mathbf G}}(s)$ is connected, as well as  $\eta(C_{{\mathbf G}}(s))$.

Observe first that $  \eta(C_{{\mathbf G}}(s))$ has finite index in  $C_{O({\mathbf V})}(\eta(s))$.  Indeed, let $M=\{g\in {\mathbf G}: \, [g,s]\in {\rm ker}\, \eta\}$, which coincides with $\{g\in {\mathbf G}: \, [\eta(g),\eta(s)]=1\}=\eta\up (C_{SO(\mathbf{V})}(\eta(s)))$. Then $\eta(M)=C_{SO(\mathbf{ V})}(\eta(s))$. As ${\rm ker}\, \eta\subseteq Z({\mathbf G})$, it follows that the mapping $g\ra [g,s]$ $(g\in  M$) is a \ho $M\ra Z({\mathbf G})$ whose kernel
is~$C_{{\mathbf G}}(s)$. The  group $Z({\mathbf G})$ is finite, so $C_{{\mathbf G}}(s)$ has finite index in $M$. So $\eta(C_{{\mathbf G}}(s))$ has finite index in $\eta(M)=C_{SO(\mathbf {V})}(\eta(s))$,  and hence in $C_{O(\mathbf {V})}(\eta(s))$.

Choose a basis $B$, say, of ${\mathbf V}$ such that $B\cap {\mathbf W}_i$ is a basis of ${\mathbf W}_i$ for $i=1,2,3$. Then, under this basis,  the matrix $t$ of $\eta(s)$ on ${\mathbf V}$ is $\diag(\Id,-\Id,s')$. Therefore,
$C_{O({\mathbf V})}(t)\subset O(\mathbf{W}_1)\times O(\mathbf{W}_2)\times C_{O(\mathbf{W}_3)}(s')$. Note that $s'\in SO({\mathbf W}_3)$ as $\dim {\mathbf W}_2$ is even (Lemma~\ref{ei2}).

As $\pm 1$ are not eigenvalues of $s'$, the group
 $C_{O(\mathbf{W}_3)}(s')$ is connected (Lemma~\ref{gf4}).
In addition, $SO(\mathbf{W}_1)\times SO(\mathbf{W}_2)\times C_{O(\mathbf{W}_3)}(s')$ is connected (as so is $C_{O(\mathbf{W}_3)}(s')$)
and has finite index in $O(\mathbf{W}_1)\times O(\mathbf{W}_2)\times C_{O(\mathbf{W}_3)}(s')$.
So both $ \eta(C_{{\mathbf G}}(s))$ and $SO(\mathbf{W}_1)\times SO(\mathbf{W}_2)\times C_{SO(\mathbf{W}_3)}(s')$ are connected subgroups of finite index in $O(\mathbf{W}_1)\times O(\mathbf{W}_2)\times C_{O(\mathbf{W}_3)}(s')$. As the connected component of the identity in an \ag is unique, these groups coincide,
as stated. $\Box$

\med
Lemma~\ref{55n} implies the following result on unipotent characters which is essential in what follows:

\bl{st5} Let $G=Spin(V)$ or $Sp(V)$, where $V$ is an orthogonal or symplectic space over $\FF_q$. Let $s\in G$ be a semisimple  element and $W_1,W_2$ be the $1$- and $-1$-eigenspaces of s on $V$.  Let $W_3=(W_1+W_2)^\perp$. Then $\nu(C_{G}(s))=\nu(SO(W_1))\cdot \nu(SO(W_2))\cdot \nu(C_{SO(W_3)}(s'))$, where $s'$ is the restriction of $s$ to $W_3$.\el

Proof. We omit the proof for $Sp(V)$ as it is straightforward. Let $G=Spin(V)$. Note that $\nu(C_{SO(W_3)}(s'))$ is meaningful as $C_{SO(W_3)}(s')$ is a finite reductive group (Lemma~\ref{f14}).
We use the notation of Lemma~\ref{55n}, assuming that $\mathbf{V}=V\otimes \overline{\FF}_q$ and that the structure of an orthogonal space on $\mathbf{V}$ is defined by the same Gram matrix as that of $V$. Then $\mathbf{W}_i=W_i\otimes \overline{\FF}_q$ for $i=1,2,3.$
Let $F$ be the Frobenius endomorphism of $\mathbf{G}$ such that $\mathbf{G}^F=G$; we keep  $F$ for the Frobenius endomorphisms of $SO(\mathbf{V})$, $SO(\mathbf{W}_i)$ $(i=1,2,3)$ inherited from
that of $\mathbf{G}$. By Lemma~\ref{55n}, $\eta(C_{{\mathbf G}}(s))=SO(\mathbf{W}_1)\times SO(\mathbf{W}_2)\times C_{SO(\mathbf{W}_3)}(s')$.
By Lemma~\ref{f11},  $\nu((\eta(C_{{\mathbf G}}(s)))^F)=\nu(C_{{\mathbf G}}(s)^F)$, and the left hand side is equal to $\nu(SO(\mathbf{W}_1)^F\times SO(\mathbf{W}_2)^F\times C_{SO(\mathbf{W}_3)}(s')^F)
=\nu(SO(W_1))\cdot \nu(SO(W_2))\cdot \nu(C_{SO(W_3)}(s'))$, as claimed. $\Box$

\section{Proof of Theorem~\ref{th1}}

 Here $G^*\cong GL_n(q)$ or $U_n(q)$. To simplify notation, we deal below with $G$ in place of $G^*$, that is, we choose a semisimple element $s\in G$ and show that the number of unipotent characters in $C_G(s)$ does not exceed $\beta(n)$.

\medskip
For our purpose, we quote the \f well known result, see \cite[p. 465]{Ca}.

\bl{gu1} Let ${\mathbf  G}=GL_n(\overline{\FF}_q)$ and $G={\mathbf  G}^F\cong GL_n(q)$ or $U_n(q)$ (depending on $F$). Then the number of unipotent characters of $G$ equals $p(n)$, the number of
partitions of $n$.\el

Let $G=GL_n(q)$, $V$ the natural $\FF_qG$-module  and $s\in G$ be a semisimple element.
We can write $V=\oplus V_i$, where
$V_i$ are the homogeneous components for~$s$, that is, each $V_i$ is a sum of isomorphic $\FF_q\lan s\ran$-modules, and distinct $V_i,V_j$ have no common \ir constituents. Let $s_i\in GL(V_i)$ be the restriction of $s$ to~$V_i$. Then $C_G(s)\subset \Pi_i GL(V_i)$, and $C_G(s)=
\Pi_i C_{GL(V_i)}(s_i)$. Let $d_i$ be the dimension of a minimal $\FF_q\lan s\ran$-submodule of $V_i$. Then $C_{GL(V_i)}(s_i)\cong GL_{d_i}(q^{m_i})$, where $m_i=\dim V_i/d_i$. One observes that the decomposition $V=\oplus V_i$ is unique up to reordering the terms. Let
$k$ be the number of terms and $n_i=\dim V_i$. Then $s$ determines the string $(n_1\ld n_k)$
up to reordering of the $n_1\ld n_k$, which is a partition of $n$, and we denote by $\pi(s)$ the partition $(n_1\ld n_k)$. (We can assume $n_1\geq \cdots \geq n_k$ but we prefer  to allow any ordering.)
If $s\in U_n(q)\subset GL_n(q^2)$ then $\pi(s)$ {\it is defined as the partition obtained for $s$ in}  $GL_n(q^2)$. The \f lemma is well known.

\bl{g22} Let ${\mathbf  G}=GL_n(\overline{\FF}_q)$, $G={\mathbf  G}^F\cong GL_n(q)$,
and let $s\in G$ be a semisimple element. Then $C_G(s)$ is isomorphic to the direct product of groups
$GL_{d_i}(q^{m_i})$, where $\sum_i d_im_i=n$.\el

The \f lemma is also well known, but we give a proof for the reader's convenience and in order to make further discussions more transparent.

\bl{u22} Let ${\mathbf  G}=GL_n(\overline{\FF}_q)$, $G={\mathbf  G}^F\cong U_n(q)$,
and let $s\in G$ be a semisimple element. Then $C_G(s)$ is isomorphic to the direct product of groups $GL_{d_i}(q^{2m_i})$ and $U_{e_j}(q^{l_j})$, where $l_i$ is odd and $\sum_i 2d_im_i+\sum e_jl_j=n$.\el

Proof. Note that each of the sums $\sum_i 2d_im_i$, $\sum e_jl_j$ can be absent. It is well known that there is an orthogonal decomposition $V=(\oplus V_i)\oplus (\oplus V_j)$, where each $V_j$ is a non-degenerate homogeneous component for $s$, and each $V_i$ is  the sum of two totally isotropic   homogeneous components for $s$. Let $H$ be the stabilizer in $G$ of this decomposition, that is,
$H=\{g\in G:gV_i=V_i, gV_j=V_j, \text{ for each term } V_i,V_j \}$.
Let $n_i=\dim V_i$, $n_j=\dim V_j$ and let $H_i,H_j$ be the restriction of $H$ to $V_i,V_j$, respectively. Then $H_j\cong U_{n_j}(q)$  and  $H_i\cong GL_{n_i/2}(q^2)$. Therefore, $n=\sum n_i+\sum n_j$ and $H\cong (\Pi_iGL_{n_i/2}(q^2))\times (\Pi _jU_{n_j}(q))$. Let $s_i,s_j$ be the restriction of $s$ to $V_i,V_j$, respectively. Then $s_i\in H_i$, $s_j\in H_j$ and
$C_G(s)=(\Pi_i C_{H_i}(s_i)) \times (\Pi_j C_{H_j}(s_j))$. Using the isomorphism  $H_i\cong GL_{n_i/2}(q^2)$ we can view a homogeneous component $V_i'$ of $V_i$ as a natural $\FF_{q^2}GL_{n_i/2}(q^2)$-module, and then $s_i$ is a homogeneous element of $GL_{n_i/2}(q^2)$, that is, $V_i'$ is a homogeneous $\FF_q\lan s_i\ran$-module.
 As in Lemma~\ref{g22}, $C_{H_i}(s_i)\cong GL_{d_i}(q^{2m_i})$, where $d_im_i=n_i/2$.
 It is also known that  $C_{H_j}(s_j)\cong U_{e_j}(q^{2l_j})$, where $e_jl_j=n_j$ and $l_j$ is odd.  So the result follows. $\Box$

\begin{lemma}\label{gu2}  Let $G=GL_n(q)$ or $U_n(q)$, and let $s\in G^*$ be a semisimple element. Then  $|\mathcal{E}_s|=\nu(C_G(s))\leq  \beta(n)$.

Furthermore, suppose that equality holds. Then $\pi(s)=\pi(n)$,
where $\pi(n)$ is defined in Theorem $\ref{th1}$, and if $G=GL_n(q)$  then $|s|$ divides $q-1$,  if $G=U_n(q)$  then $|s|$ divides $q+1 $.
\end{lemma}

Proof.  If $G=GL_n(q)$ then, by Lemma~\ref{g22},  $C_G(s)\cong \Pi_i GL_{d_i}(q^{m_i})$, where $\sum_i d_im_i=n$. Recall (Lemma~\ref{gu1}) that the number of unipotent characters of $GL_n(q)$ equals $p(n)$ and hence does not depend on $q$. So  $\nu(C_{G^*}(s))=\Pi_i p(d_i)$. Set $n'=\sum d_i$. Then $\Pi_i p(d_i)\leq \beta(n')$. By Lemma~\ref{bo1},
$\beta(n')<\beta(n)$ for $n'<n$; if equality holds above, then $n=n'$, and hence $d_i=n_i$ for every $i$.
This implies that each $V_i$ is a sum of one-dimensional $s$-stable subspaces, whence the result.

Let $G=U_n(q).$ Then $C_{G}(s)$ is a direct product
of groups isomorphic to $GL_{m_i}(q^{2d_i})$,  $i=1,\ldots ,k'$, and $U_{l_j}(q^{f_j})$,  $j=1,\ldots ,k''$,
for some integers $k',k''\geq 0$, and  $n=\dim V= 2\sum_{i=1}^{k'}m_id_i+\sum_{j=1}^{k''} e_jf_j$. (Note that
$C_{G}(s)$ may be a product of $GL_{m_i}(q^{2d_i})$ or  $U_{l_j}(q^{f_j})$ only.)
The number of unipotent characters of $GL_{m_i}(q^{2d_i})$  equals $p(m_i)$ and that of $U_{l_j}(q^{f_j})$   equals $p(l_j)$ (Lemma~\ref{gu1}). Let $n'=\sum m_i$, $n''=\sum l_j$.
Then $|\mathcal{E}_s|=\nu(C_G(s))=\prod p(m_i)\cdot\prod p(l_j)\leq \beta(n')\cdot\beta(n'')$.
By Lemma~\ref{bo1}, $\beta(n')\cdot\beta(n'')\leq \beta(n'+n'')$. If the equality holds then $n=n'+n''$, whence $n'=0$, $n=n''$ and $f_j=1$ for all $j=1\ld k''$. \itf $|s|$ divides $q+1$ and $(l_1\ld l_{k''})=\pi(n)$, so  the result follows as above. $\Box$

\med
We now show that the bound is attained for every $n$ for $q$ large enough.

\begin{lemma}\label{gu21}
Let $n,i\in\NN$, $i\in \{0,1,2,3\}$ with $i\equiv n\pmod 4$. Assume that $n>3$ if $i=0,1,2$, and $n>10$ for $i=3$.

Let $G=GL_n(q)$, respectively, $U_n(q)$. If  $n\leq 4(q-1)+i$, respectively, $n\leq 4(q+1)+i$, then  $\nu(C_G(s))= \beta(n)$ for a suitable semisimple element $s\in G$.
\end{lemma}

Proof. Let $n=4k+i$. Then $k\leq q-1$, respectively, $q+1$. Therefore, there exist distinct  elements  $a_1\ld a_k\in GL_1(q)$, respectively, $U_1(q)$. If $i=0$ then we
 set $s=\diag(a_1\cdot \Id_4,a_2\cdot \Id_4\ld  a_k\cdot\Id_4 )$.
If $i=1$  then we take the last scalar to be $a_k\cdot\Id_5$, if $i=2$ then we take the last scalar to be $a_k\cdot\Id_6$. If $i=3$ then  we take the last two scalars to be $a_{k-1}\cdot\Id_5$ and $a_k\cdot\Id_6$.   If $G=U_n(q)$ then we choose  an orthogonal basis of the underlying space, in order to get $s\in U_n(q)$. Then $C_{G}(s)$ is the direct product of groups $GL_4(q)$ (respectively $U_4(q)$) if $n\equiv 0\pmod 4$, with obvious adjustments in the other cases. Then  $\nu(C_{G^*}(s))=\beta(n)$. So the bound $\beta(n) $ is attained. $\Box$

\medskip

{\it Proof of Theorem} \ref{th1}. This follows from Lemmas \ref{gu21} and \ref{gu2}. $\Box$

\bl{cg1} Let C be a cyclic group, $|C|>2$. Set $l=(|C|-2)/2$ if $|C|$ is even, and $l=(|C|-1)/2$ if $|C|$ is odd. Then there are l distinct  elements  $a_1\ld a_l\in C$ such that $a_ia_j\neq 1$ for all  $1\leq i,j \leq l$. \el

Proof. Let $C=\lan a \ran$. Then set $a_i=a^i$. As the  elements  $a^i$ $(1\leq i\leq |C|-1)$ are all distinct, and $a^{l+1}$ is of order 2 if $|C|$ is even, it follows that $\{a^j: 1\leq j\leq l\}$ satisfies the conclusion of the lemma. $\Box$

\med
For application to other classical groups we need a slightly different version of Lemma~\ref{gu21}.
We view $GL_n(q)$ as a matrix group over $\FF_q$ and  $U_n(q)$ as a matrix group over $\FF_{q^2}$ whose subgroup of diagonal matrices is $\diag(U_1(q)\ld U_1(q))$. In Lemma~\ref{cg2} below $D$ denotes
the group of diagonal matrices in $G$. For $d\in D$ the set of distinct diagonal entries of $d$ is denoted by $Spec(d)$.

\bl{cg2} Let $G=GL_n(q)$ or $U_n(q)$  and let $G_2$ be  the subgroup of  $G$ of index~$2$ if $q$ is odd, and $G_2=G$ if $q$ is even. Suppose that $q\geq n+5$. Then there exists a semisimple element
$s\in D\cap G_2$ such that $Spec(s)\cap Spec(s\up)=\emptyset$ and $\nu(C_G(s))=\beta(n)$. \el

Proof.  Let $C=GL_1(q)$ or $U_1(q)$ if $q$ is even, and let $C$  be the subgroup  of index~2  in these groups if $q$ is odd. Let $l$ be as in Lemma~\ref{cg1}. Then $l=(q-2)/2$ if $q$ is even,
$(q-3)/2$ if $q\equiv 3\pmod 4$ and $(q-5)/2$ if $q\equiv 1\pmod 4$.  By  Lemma~\ref{cg1}, for every $k\leq l$ there are  distinct  elements  $a_1\ld a_k\in C$ such that $a_ia_j\neq 1$ for  all
$1\leq i,j \leq k$.

Then we take $k=(n-r)/4$, where $0\leq r<4$ and $n\equiv r\pmod 4$.
As $q\geq n+5$ by assumption, we have $k=(n-r)/4 \le (q-5)/4 \le l$.
Let us choose  these  elements $a_1\ld a_k$ for a similar reasoning as in the proof of
Lemma~\ref{gu21}
to construct suitable elements $s\in D$.

   Then   $\nu(C_G(s))=\beta(n)$ by Lemma~\ref{gu1}.
 In addition,  as $s$ is a diagonal matrix with entries $a_1\ld a_k$ (with certain multiplicities) the condition $a_ia_j\neq 1$ for all
 $1\leq i,j\leq k$
 implies $Spec(s)\cap Spec(s\up)=\emptyset$. As each diagonal entry of $s$ lies in~$C$, it follows that  $s\in G_2$.  $\Box$

\section{Other classical groups}

\subsection{Remarks on classical groups}

We start with observations on  the centralizers of semisimple elements of classical groups. Let $H\in\{O_{2n+1}(q), $, $q$ odd, $Sp_{2n}(q), O^{\pm}_{2n}(q)\}$ and let $V$ be the underlying space for $H$. Recall that $\Omega^\pm_{2n}(q)$ denotes the subgroup of  $O^{\pm}_{2n}(q)$ formed by elements of spinor norm 1, and in even characteristic  $\Omega^\pm_{2n}(q)=SO^\pm_{2n}(q)$ by convention.

\medskip

The following two lemmas are well known.

\begin{lemma}\label{kk2}   \cite[Prop. 2.5.13]{KL} For $q$ odd, set $\ep(n)=(-1)^{(q-1)n/2}$. The group  $\Omega^+_{2n}(q)$, respectively,  $\Omega^-_{2n}(q)$ contains
$-\Id$ \ii $\ep(n)=1$, respectively, $\ep(n)=-1$. In particular, $\Omega^+_{2n}(q)$  contains
$-\Id$ if $n$ is even or $q$ is a square. \end{lemma}

\begin{lemma}\label{ei2}  Let  $G\in \{SO_{2n+1}(q), q ~odd, SO^\pm_{2n}(q),~Sp_{2n}(q)\},$   and let $V$ be the natural module for G. Let $g\in G$ be a semisimple element and let $V_1$ and $V_2$  be the   $1$- and $-1$-eigenspaces of g on V. (If q is even then $V_2=0$ by convention.) Then

$(1)$ $V_1$ and $V_2$  are non-degenerate and orthogonal to each other;

$(2)$ $\dim V_2$ and $\dim(V_1+V_2)^\perp$ are even;

$(3)$ $\dim V_1$ is even unless $G= SO_{2n+1}(q), $ in which case $\dim V_1$ is odd.
\end{lemma}

Proof.  Let $i\in\{1,2\}$. (1) If $V_i$ is degenerate then $U:=V_i\cap  V_i^\perp\neq 0$ is totally isotropic.
Let $0\neq u\in U$, so $\dim u^\perp=\dim V-1$ \cite[2.1.5]{KL}.  As $g$ is semisimple, $u^\perp$ has a $g$-invariant complement $U', $ say. Let $v\in U'$ and let
$f$ be the form on $V$ defining $G$.
Then $0\neq f(u,v)=f(gu,gv)=af(u,gv)$, where $a=1$ or $-1$. \itf $gv=av$, which is a contradiction as such a $v$ must be in $V_i$.

If $V_1,V_2\neq 0$ then $q$ is odd; choose $0\neq v_i\in V_i;$ then $f(v_1,v_2)=f(gv_1,gv_2)=-f(v_1,v_2)$, whence $f(v_1,v_2)=0$.

(2) It suffices to prove this statement for the respective groups over $\overline{\FF}_q$; in this case
$V$ is the sum of the eigenspaces of $g$, and $\pm 1$ are not \eis of $g$ on $W:=(V_1+V_2)^\perp$.
Let $e$ be an \ei of $g$ on $W$, so $e\neq \pm 1$, and  $W_e$ be the respective eigenspace. Then for $0\neq w\in W_e$
we have $f(w,w)=f(gw,gw)=e^2f(w,w)=0$ as $e^2\neq 1$. One easily observes that $w^\perp=\lan w\ran +W'$, where $W'$ is a $g$-stable non-degenerate subspace of  $w^\perp$. By induction,
$\dim W'$ is even, and hence so is $\dim W$.

Moreover, if $v\notin w^\perp$
then, as in the proof of (1),
$gv=e\up v +x$
for $x\in w^\perp$. This implies by induction that the determinant of $g_W$, the restriction of $g$ to $W$,
equals~1.  As $\det g=1$ and $g$ acts on $V_2$ as $-\Id$, it follows that $\dim V_2$ is even, as claimed.

(3) is obvious as  $\dim V_1=\dim V-\dim V_2-\dim W$. $\Box$

\med
Next we describe the structure of centralizers of semisimple elements in $H$. This is treated in \cite[\S 1]{FS} and elsewhere, but we choose to briefly recall the main facts in a form compatible with what follows.

Let $h\in H$ be a semisimple element. Viewing $V$ as $\lan h\ran$-space we can write $V=V_1\oplus\cdots \oplus V_k\oplus V_{k+1}\oplus\cdots \oplus V_{k+l}$, where $V_1\ld V_{k+l}$
are homogeneous components of $V$ for $\lan h\ran$. (In other words, $V_1\ld V_{k+l}$
 are $h$-stable, for every $i\in \{1\ld k+l\}$ all \ir constituents of $V_i|_{\lan h\ran}$
are isomorphic to each other and not isomorphic to those of $V_j|_{\lan h\ran}$ for every $j\neq i$.) Furthermore, each $V_i$ is either non-degenerate or totally isotropic, see for instance \cite[Lemma 3.3]{SZ}. By reordering the terms, we assume  that $V_1\ld V_k$ are totally isotropic (unless $k=0$) whereas $V_{k+1},\ldots ,V_{k+l}$ are non-degenerate (unless $l=0$). In the former case for every $i\leq k$ there is another   totally isotropic homogeneous component
$V_j$, say, such that $V_i|_{\lan h\ran}$ and $V_j|_{\lan h\ran}$ are dual to each other and  $V_i+V_j$ is non-degenerate  \cite[Lemma 3.3]{SZ}. It follows that $k=2m$ is even.  We can reorder $V_1\ld V_{k}$ so that $V_{i},V_{k-i+1}$ are dual as $\lan h\ran$-modules, $i=1\ld m$. Set
$h_i=h|_{V_i}$ and $n_i=\dim V_i$ for $i=1\ld k+l$. If $h_i=\pm \Id$ then $V_i$ is non-degenerate (see Lemma~\ref{ei2}), and hence $i>k$ in this case.

For  $i\in \{1\ld k\}$ set $H_i=GL(V_i)$ and for $i\in \{k+1\ld l\}$ set $H_i=\{g\in H: g|_{V_j}=\Id$ whenever $j\neq  i\}\cong I(V_i)$. (For uniformity, we use $I(V_i)$ to denote the classical groups defined by the relevant form on $V_i$.)
 Then $$C_H(h)\cong C_{H_1}(h_1)\times \cdots \times C_{H_m}(h_m)\times C_{H_{k+1}}(h_{k+1})\times \cdots \times C_{H_{k+l}}(h_{k+l}).$$

Let  $d_i$ be the dimension of each \ir constituent of $h_i$, $i=1\ld k+l$. As $V_i$ is homogeneous,   $n_i$ is a multiple of $d_i$. Write $n_i=d_ie_i$. If $i\leq k$ then $C_{H_i}(h_i)=C_{GL(V_i)}(h_i)$.

\medskip
(a) Suppose that $H$ is symplectic. If $h_i=\pm \Id$ then $C_{H_{i}}(h_{i})
\cong Sp_{n_i}(q)$. If  $i\leq k$ then $C_{H_i}(h_i)\cong GL_{e_i}(q^{d_i})$;
 if $i>k$ and $h_i\neq \pm \Id$  then $C_{H_i}(h_i)\cong U_{e_i}(q^{d_i/2})$. (Here we write
$U_{e_i}(q^{d_i/2})$ due to our notation for unitary groups, that is, $U_{e_i}(q^{d_i/2})\subset GL_{e_i}(q^{d_i})$.)

\medskip
(b) Suppose that $H$ is orthogonal.  If $h_i=\pm \Id$ then $C_{H_{i}}(h_{i})=H_i
\cong O(V_i)$. If $h_i\neq \pm \Id$  and $i\leq k$ then $C_{H_i}(h_i)\cong GL_{e_i}(q^{d_i})$. If $h_i\neq \pm \Id$  and $i>k$   then $C_{H_i}(h_i)\cong U_{e_i}(q^{d_i/2})$, where $e_i$ is odd
\ii the Witt defect of $V_i$ is 1.

In case (b)  fix some $V_i$ of Witt defect 1 (assuming the existence of it). Then $V_i$ is a direct sum of $e_i$  \ir non-degenerate $\langle  h_i\rangle$-modules isomorphic to each
other.
Denote by $D$ one of them, so $h_i$ acts irreducibly  on $D$. Here $\dim D>1$ as $h_i\neq \pm 1$. Therefore the Witt defect of $D$ is $1$ because
otherwise  $O(D)$ has no \ir  element (\cite[Satz 3(c)]{Hu}). So the assertion on the parity of  $e_i$ follows from  \cite[Proposition 2.5.11(ii)]{KL}.

(Note that $d_i=\dim D/2$ can be even.)

\medskip
We  state the above information in a uniform way as follows:

\begin{propo}\label{jj6} Let $h\in H$ be a semisimple element and let $V_1,V_2$ be the $1$- and $-1$-eigenspace of $h$ on $V$. Then $C_H(h)\cong I(V_1)\times I(V_2)\times \Pi_i GL_{d_i}(q^{l_i})\times \Pi_j U_{e_j}(q^{m_j})$, where $\frac{1}{2}(\dim V_1+\dim V_2)+\sum _id_il_i+\sum e_jm_j=n$.
\end{propo}

\begin{corol}\label{jj7} Let $G\in \{SO_{2n+1}(q), q ~odd, SO^\pm_{2n}(q),~Sp_{2n}(q)\},$   and let $V$ be the natural module for G. Let $s\in G$ be a semisimple element. Suppose that s does not have \eis $-1$ on V and the \mult of the \ei $1$ is at most $ 1$. Then $C_G(s)\cong\Pi_i GL_{d_i}(q^{l_i})\times \Pi_j U_{e_j}(q^{m_j})$, where $ \sum _id_il_i+\sum e_jm_j=n$.
\end{corol}

Proof. Suppose that $\dim V$ is even. Then, under these assumptions, $C_G(s)=C_H(s)$ by Lemma~\ref{f14}(2), so the result follows from Proposition \ref{jj6}. If $\dim V$ is odd then $O(V)=SO(V)\times \{\pm \Id\}$, so $C_H(s)=C_G(s)\times \{\pm \Id\}$.
 $\Box$

\begin{lemma}\label{ed1}  Let $s\in G=SO_{2n}^-(q)$  be a homogeneous semisimple element, and $s\neq \pm \Id$. Then $C_G(s)\cong U_e(q^d)$,  where $ed=n$, e is odd,  and   $\nu(C_G(s))\leq p(n')$, where $n'$  is the greatest odd divisor of n. In addition, if $(n',q)\neq (n,3)$ then there exists a  (homogeneous) semisimple element $s'\in \Omega_{2n}^-(q)$ such that $C_{G}(s')\cong U_{n'}(q^{n/n'})$.\end{lemma}

Proof.  By the comment prior to Proposition \ref{jj6} and Lemma~\ref{f14}(2), we have
$C_G(s)\cong U_e(q^d)$, where $e$ is odd and $n=de$. By Lemma~\ref{gu1},
$\nu(U_e(q^d))=p(e)$, and $p(e)\leq p(n')$. For the additional claim, decompose
the natural $\FF_qG$-module $V$ as a direct sum of $n'$ non-degenerate subspaces of dimension $2n/n'$ and of Witt defect 1. Let $D$ be one of them. Then $SO(D)\cong SO_{2n/n'}^-(q)$, so $SO(D)$ contains an \ir element $t$, say, of order $q^{n/n'}+1$ \cite{Hu}.  Then $t^2$ is still \ir on~$D$ unless $n=n'$ and $q= 3$. Choose $s$ to be an element of $G$ stabilizing each direct summand (which is isomorphic to $D$) and
acting on each of them  as $t^2$ does. Then $s$ is homogeneous and $C_G(s)\cong U_{n'}(q^{n/n'})$ by the above. So the claim follows. $\Box$

\subsection{Subgroups of classical groups and their unipotent characters}

We assume the group  ${\mathbf G}^*$ to be simply connected, which in turn
guarantees $C_{{\mathbf G}^*}(s)$ to be connected for every semisimple
element $s\in {\mathbf G}^*$, see Lemma~\ref{gf4}(2). In view of Lemma~\ref{f13} our task is to  obtain a sharp upper bound for $\nu(C_{G^*}(s))$.
The information on the number of unipotent characters of $G$ is given in \cite[Section~13.8]{Ca}.

\medskip
If $G^*=Spin_{2n+1}(q)$, $q$ odd, or $Spin^\pm_{2n}(q)$  then the natural module $V$, say,
for $O_{2n+1}(q)$ or $O^\pm_{2n}(q)$ can be viewed as $\FF_qG^*$-module under
the natural \ho of $G^*$ into the respective classical group.
So we refer to $V$ as the natural module for $G^*$.

\medskip
The function $\beta(n)$ plays a significant role in this paper.  It is not true that $|{\mathcal  E}_s|\leq \beta(n)$, but the \f lemma singles out an important special case where this is true.

\begin{lemma}\label{ei1}
 Let  $G^*\in \{ Spin_{2n+1}(q)$ for $q$ odd, $ Spin^\pm_{2n}(q)$,   $Sp_{2n}(q)\}, $ and let $V$ be the natural module for $G^*$. Let $s\in G^*$ be  a  semisimple element   such that
  the \mult of \eis   $1$  and  $-1$ of $s$   on $V$  does not exceed $1$.

$(1)$   $|{\mathcal E}_s|=\nu(C_{G^*}(s))\leq \beta(n)$.

$(2)$ If $~V=V'\oplus V''$ is an orthogonal decomposition, such that $sV'=V',sV''=V''$ and $Hom_s(V',V'')=0$
(equivalently, $s$   has no common \ei on ${\mathbf V}',{\mathbf V}''$) then
  $\nu(C_{G^*}(s))=\nu(C_{SO(V')}(s_1))\cdot \nu(C_{SO(V'')}(s_2))$,
where $s_1,s_2$ are the restriction of $s$ to $V',V''$, respectively.
\end{lemma}

Proof.  By Lemma~\ref{ei2},  the \mult of the \ei $-1$ is always even,  as well as of the \ei 1 unless $G^*\cong Spin_{2n+1}(q)$, where the \mult of the \ei 1 is always odd. Therefore, the assumption implies that $-1$ is not an \ei of~$s$, as well as 1, provided $G^*\neq Spin_{2n+1}(q)$. By  Lemma~\ref{gf4}, $C_{{\mathbf  G}^*}(s)$ is connected; so by Lemma~\ref{f13}, $|{\mathcal E}_s|=\nu(C_{G^*}(s))$.

(1) Let ${\mathbf  H}\in \{SO_{2n+1}(\overline{\FF}_q)$,  $q$ odd, $SO_{2n}(\overline{\FF}_q)$, 
$Sp_{2n}(\overline{\FF}_q)\}$, and $\eta:{\mathbf  G}^*\ra {\mathbf  H}$ the natural homomorphism. Keep $F$ to denote the Frobenius endomorphism of ${\mathbf  H}$ inherited from that of ${\mathbf  G}^*$, and set $H={\mathbf  H}^F$. Then $\eta$ is surjective and
 $H$ is one of the groups $SO_{2n+1}(q)$, $ q$ odd, $  SO^\pm_{2n}(q)$, 
 $Sp_{2n}(q) $ (depending on $G^*$ and $F$).
As $-1$ is not an \ei of $\eta(s)$, by Lemma~\ref{gf4}(3), the group $C_{{\mathbf H}}(\eta(s))$ is connected and, by Lemma~\ref{f14},
$\nu(C_{ G^*}(s))=\nu(C_{H}(\eta(s) ))$. By Corollary \ref{jj7}, we have $C_{G^*}(s)\cong
\Pi_i GL_{d_i}(q^{l_i})\times \Pi_j U_{e_j}(q^{m_j})$, where $ \sum _id_il_i+\sum e_jm_j=n$. Furthermore,
the number of unipotent characters of each factor is equal to $p(d_i)$ or $p(e_j)$
(Lemma~\ref{gu1}), so the total is $\Pi_i p(d_i)\cdot \Pi_j p(e_j)$.
By~\cite{BO},  this number is not greater than $\beta(n)$, whence (1).

(2) By Lemma~\ref{gf4}(3), the group $C_{O(\mathbf{V}_i)}(s_i)=C_{SO(\mathbf{V}_i)}(s_i)$ is  connected for $i=1,2$, so $C_{O(\mathbf{V})}(s)=C_{SO(\mathbf{V})}(s)=C_{SO(\mathbf{V}_1)}(s_1)\times C_{SO(\mathbf{V}_2)}(s_2)$. In addition, $C_{SO(\mathbf{V}_i)}(s_i)^F=C_{SO(V_i)}(s_i) $, so $ C_{SO( V)}(s)=C_{SO( V_1)}(s_1)\times C_{SO( V_2)}(s_2)$. This implies (2).
 $\Box$

\medskip
The \f lemma tells us that the bound in Lemma \ref{ei1}(1) is attained if $q$ is large enough and
$G^*\neq SO^-_{2n}(q)$.

\bl{rr1} Let $G^*\in\{SO_{2n+1}(q)$, $q$ odd,
$Sp_{2n}(q)$, $SO^+_{2n}(q)\}$, and let $V$ be the natural $\FF_qG^*$-module.
Suppose that $n\leq  q-5$. Then there exists $t\in G^*$ such that $V$ is the sum of the eigenspaces of $t$, the \mult of the \eis $1$ and $-1$ of $t$ is at most $1$ and $\nu(C_{G^*}(t))=\be(n)$.
In addition, if $q$ is odd and $G^*$ is orthogonal 
then $t$ can be chosen in a subgroup of index $2$ of $G^*$. \el

Proof. It is well known that
there exist totally singular subspaces $V_1,V_2$ of $V$ such that  $V_1\cap V_2=0$,
$\dim V_1=\dim V_2=n$, $V_1+V_2$ is non-degenerate, and there are dual bases in $V_1,V_2$
in the sense that if $g\in G$ with $gV_i=V_i$ $(i=1,2)$ and $g_i$ is the matrix of $g$ on $V_i$, then $g_2={}^Tg_1\up$, where ${}^Tg_1$ is the transpose of $g_1$. Moreover, for $H=GL(V_1)\cong GL_n(q)$ there is an embedding $\lam:H\ra G$ such that $\lam(h)=\diag(h,{}^Th\up)$ or $\diag(h,1,{}^Th\up)$ for $h\in H$. Let $W$ be the natural module for $H$. \itf $V_1|_H\cong W$ and $V_2|_H$ is dual to $W$.

Let $s\in H$ be as in Lemma~\ref{cg2} and $t=\lam(s)$. Then the statement on the \eis of $t$ on $V$ is obvious. Since $W$ is the sum of the eigenspaces of $s$, it follows that $V$ is the sum of the eigenspaces of $t$. In addition, the choice of $s$  in Lemma~\ref{cg2} implies  every eigenspace of $\lam(s)$ to
lie in $V_1$ or $V_2$. It easily follows that $C_{G^*}(t)=\lam(C_H(s))\cong C_H(s)$. Therefore,
the number of unipotent characters of $C_{G^*}(t)$ and $ C_H(s)$ is the same. By Lemma~\ref{cg2}, the latter is equal to $\be(n)$, whence the result.

However, to be precise, the isomorphism $C_H(s)\ra C_{G^*}(t)$ should be accompanied with an isomorphism of algebraic groups $C_{{\mathbf H}}(s)\ra C_{{\mathbf G}^*}(t)$ such  that
$C_{{\mathbf H}}(s)^F=C_H(s)$ and  $ C_{{\mathbf G}^*}(t)^F=C_{G^*}(t)$.   (As above, we use the same letter $F$ for the Frobenius endomorphism of different groups $C_{{\mathbf H}}(s)$ and $ C_{{\mathbf G}^*}(t)$).

Let ${\mathbf G}^*= SO_{2n+1}(\overline{\FF}_q)$, $ Sp_{2n}(\overline{\FF}_q)$ or $SO_{2n}(\overline{\FF}_q)$ and ${\mathbf H}=GL_n(\overline{\FF}_q)$. In each case
we choose for $F$ the standard   Frobenius endomorphism arising from
raising matrix entries of elements of the above groups to the $q$-power (see \cite[p. 37]{DM}). (For this we choose a basis $B$ in $V$ as above, such that $B\cap V_1$ and $B\cap V_2$ are dual bases, and view it as a basis of  the underlying space ${\mathbf V}$ of ${\mathbf G}$.) Then
  $G^*={\mathbf G}^{*F}$ and $H={\mathbf H}^F$. The latter holds true when we consider ${\mathbf H}$ as $GL({\mathbf V}_1)$ or as a subgroup of ${\mathbf G}^*$ stabilizing ${\mathbf V}_1$ and ${\mathbf V}_2$.  Then  $C_{{\mathbf H}}(s)$ and    $C_{{\mathbf G}^*}(t)$  are isomorphic, as the \eis of $s$ on ${\mathbf V}_2$ are  the inverses of those on   ${\mathbf V}_1$ (see Lemma~\ref{cg2}).
In addition, we  have $C_{{\mathbf H}}(s)^F=C_H(s)$ and    $C_{{\mathbf G}^*}(t)^F=C_{G^*}(t)$.

For the additional statement for $q$ odd let $G_2$ be the subgroup of index 2 in $G^*$. Then   $|\lam(H):(\lam(H)\cap G_2)|\leq 2$. By  Lemma~\ref{cg2}, $s$ can be chosen in the subgroup of index 2 in $H$, whence the claim. $\Box$

\medskip
 Now we consider the case where  $G^*= Spin^- _{2n}(q)$. We shall see that the statement of Lemma~\ref{rr1} remains true for $n$  odd but fails otherwise. Recall that  $\beta'(n)=\max_{a~ odd } \,\beta (a)\beta (n-a)$; by Lemma~\ref{nk8},  $\beta'(n)=\beta(5)\beta(n-5)$ for $n>6$.

\begin{lemma}\label{ei2a}  Let  $G^*=Spin^-_{2n}(q)$, and let $V$ be the natural module for $G^*$. Let $s\in G^*$ be a semisimple element such that $1$ and $-1$ are not \eis of s  on V.

 $(1)$   $\nu(C_{G^*}(s))\leq \beta(n)$ for n odd,  and
 $\nu(C_{G^*}(s))\leq \beta(5)\beta(n-5)  $ for $n> 6$ even.

 $(2)$ Suppose that $q\geq n+5$. Then the bounds in $(1)$ are attained for some $s$.

 $(3)$ If $n=6,4,2$ then the maximum of  $\nu(C_{G^*}(s))$ equals $ 9,3,2$, respectively.
\end{lemma}

Proof.  Let $\eta:G^*\ra SO(V)$ be the natural \ho and $t=\eta(s)\in \Omega(V)\cong \Omega_{2n}^-(q)$. As observed in Lemma~\ref{ei1}, $\nu(C_{G^*}(s))=\nu(C_{SO(V)}( t))$.

 If $n\leq 7$ the claims follow by inspection, whence (3).
Suppose that $n>7$.

(2) Let $V=V_1\oplus V_2$, where $V_1,V_2$ are non-degenerate subspaces of $V$, the Witt defect of $V_1$ is 1, the Witt defect of $V_2$ is 0, and $\dim V_1=10$, $\dim V_2=2n-10$.
By  Lemma~\ref{rr1}, there is an element $s_2\in \Omega(V_2)\cong \Omega^+_{2n-10}(q)$ such that $V_2$ is the sum of eigenspaces of $s_2$ (whence $s_2^{q-1}=1$) and
$\nu(C_{SO(V_2)}(s_2))=\be(n-5)$.

Furthermore, there is a homogeneous element $s_1\in \Omega(V_1)\cong \Omega^-_{10}(q)$ such that $|s_1|>2$ divides $q+1$ and  $C_{SO(V_1)}(s_1)\cong U_5(q)$ (see Lemma~\ref{ed1}).
 As $s_2^{q-1}=1$, it follows that  $s_1$, $s_2$ have no common \eis over $\overline{\FF}_q$. Let $t=\diag(s_1,s_2)$ and $s\in G^*$ be such that $\eta(s)=t$ (such $s$ exists as $\eta(G^*)=\Omega(V)$).    Then  $C_{O(V)}(t)=C_{O(V_1)}(s_1)\times C_{O(V_2)}(s_2)$; it follows from Lemma~\ref{gf4}(3) that $C_{SO(\mathbf{V})}(t)=C_{SO(\mathbf{V}_1)}(s_1)\times C_{SO(\mathbf{V}_2)}(s_2)$, and
 also that $C_{SO(V)}(t)$, $C_{SO(V_1)}(s_1)$ and $ C_{SO(V_2)}(s_2)$ are finite reductive groups.
So  $\nu(C_{SO( V)}(t))=\nu(C_{SO(V_1)}(s_1))\cdot \nu(C_{SO( V_2)}(s_2))=\be(5)\be(n-5)$.
 By Lemma~\ref{st5},   $\nu(C_{G^*}(s))=\be(5)\be(n-5)$, so we are done if $n$ is even
 and $n>6$. If  $n$ is odd, then $n-5 \equiv 0$ or  $2 \mod 4$; in both cases
 $\be(5)\be(n-5)=\be(n)$ by Lemma~\ref{bo1}, provided $n-5\ge 4$,
whence the result.

(1) If $n$ is odd, this is already proven in Lemma~\ref{ei1}.
Suppose that $n$ is even. Suppose the contrary, and let $s\in G^*$ be such that  $\nu(C_{G^*}(s))>\beta(5)\beta(n-5) $.

Choose a decomposition $V=(\oplus_{i=1}^k V_i)\oplus (\oplus _{j=1}^l V_j)$ described after Lemma~\ref{ei2},
in particular, each term is a minimal non-degenerate $s$-stable subspace of $V$,  each $V_j$ $(j=1\ld l)$ is minimal and each $V_i$ $(i=1\ld k)$ is the sum of two minimal $s$-stable subspaces of $V$. By Corollary \ref{jj7}, $C_G(s)\cong
\Pi_{i=1}^k GL_{d_i}(q^{l_i})\times \Pi_{j=1}^l U_{e_j}(q^{m_j})$, where $\dim V_i=2d_il_i$, $\dim V_j=2e_jm_j$, so $ \sum _id_il_i+\sum e_jm_j=n$. Note that each $V_i$ has Witt defect 0. By \cite[2.5.11]{KL}, at least one $V_j$ has Witt defect 1, in particular, $l\neq 0$.

Observe first that the case $k=0,l=1$ does not hold. Indeed, otherwise  $n=e_jm_j$ and $C_{G^*}(s)\cong U_{e_j}(q^{m_j})$. By Lemma~\ref{ed1},   $\nu(C_{G^*}(s))\leq p(n')$, where $n'$ is the odd part of $n$.   Then $p(n')\leq p(n/2)$, as $n$ is even, and $p(n/2)\leq \beta(n/2)$.
By Lemma~\ref{bo3}, $\beta(n/2)<7\beta(n-5)$  and $\beta(n/2)<3\beta(n-3)$.
In the latter case if $n-3\equiv 1~\pmod 4$ then $ 3\beta(n-3)=3\cdot 7\cdot 5^{(n-8 )/4}$, and
this is less than $ \beta(5)\cdot \beta(n-5)=7\cdot 77\cdot 5^{(n-16 )/4}$ (as $n-5\equiv 3\pmod 4$). This is a contradiction.

Choose $j$ so that the Witt defect of $V_j$ is 1. Set $W=V_j^\perp$, so $W\neq 0$ is the sum of all terms in the above decomposition but $V_j$. Then the Witt defect of $W$ equals 0.
 Let $s_j,s'$ be the restriction of $s$ to $V_j$, $W$, respectively. We show that $n_j:=e_jm_j$ is odd. Indeed, by Lemma~\ref{ei1}, $\nu(C_{SO(V)}(s))=\nu(C_{SO(W)}(s'))\cdot \nu(C_{SO(W)}(s'))$, and $\nu(C_{SO(W)}(s'))\leq
\beta(n-n_j)$ by Lemma~\ref{rr1}.  If $n_j$ is even then, by Lemma~\ref{ed1},    $\nu(C_{SO(V_j)}(s_j))=p(n'_{j})$, where $n_j'$ is the odd part of $n_j$. By Theorem~\ref{bo2},
$p(n'_{j})\leq \beta(n_j/2)$. By Lemma~\ref{bo3}, $\beta(n_j/2)\leq 3\beta(n_j-3)$. Then
$3\beta(n_j-3)\cdot  \beta(n-n_j)\leq 3\beta(n-3)$ by Lemma~\ref{bo1}.
 By the above, this is less than $\beta(5)\beta(n-5)$.

So $n_j$ must be odd, and hence
$\nu(C_{G^*}(s))\leq \beta'(n)=\max _{a~odd}\beta(a)\be(n-a)$. If $n>6$ then $\beta'(n)=\beta(5)\be(n-5)$ by Lemma~\ref{nk8}, as required. $\Box$

\medskip
Recall that  $\al(n),\al^+(n),\al^-(n)$ denote the number of unipotent characters of the group
$Sp_{2n}(q)$, $Spin_{2n}^+(q)$, $Spin_{2n}^-(q)$, respectively. Note that
$\nu(Spin_{2n+1}(q))=\al(n)$ as well.

An essential role in what follows is played by Lemmas \ref{2u1} and \ref{2u2}   which in a sense
 generalize Lemma~\ref{gu1} to other classical groups.

\begin{lemma}\label{2u1}  Let $G^*\in \{  Spin_{2n}^+(q),Spin_{2n}^-(q), Sp_{2n}(q)\}$, where $q$ is even. Let $a,c\geq 0$ be integers such that $a+c=n$ 
  and if $G^*=Spin_{2n}^-(q)$ then $c<n$. If $q\geq n-a+5$  then there exists a semisimple element $s\in G^*$ such that

\medskip
\begin{center}
$\nu(C_{G^*}(s))=\begin{cases}\al(a)\beta(c) &\text{if } ~G^*=Sp_{2n}(q);\\
\al^+(a)\beta(c) &\text{if } ~G^*=Spin^+_{2n}(q);\\
\al^-(a)\beta(c)&\text{if } ~G^*=Spin^-_{2n}(q).
\end{cases}$
\end{center}
\end{lemma}

Proof.  Let $V$ be the natural module for   $G^*$.
(Note that $Spin(V)\cong \Omega(V)$.) 
Then $V$ contains a non-degenerate subspace $W$, say, of dimension~$2c$ and of Witt defect~$0$.
Set $H=\{g\in  G^*: gx=x$ for every $x\in W^\perp\}$. Then  $H\cong Spin^+_{2c}(q)$ or $Sp_{2c}(q)$. By Lemma~\ref{rr1},  there is an element $h\in   H$ such that  $\nu(C_H(h))=\beta(c)$ and $h$ does not have \ei~1 on~$W$.
Then $C_{G^*}(h)\cong C_{H}(h)\times X$, where $X\cong Sp_{2a}(q)$ or $ \Spin^\pm _{2a}(q)$.
By Lemma~\ref{ei1},
$\nu(C_G(h))=\nu(C_{H}(h))\cdot \nu(X)=\beta(c)\cdot  x$, where $x=\al^+(a)$,  $\al^-(a)$ or $\al(a)$ when $G^*\cong Spin_{2n}^+(q),Spin_{2n}^-(q), Sp_{2n}(q)$, respectively. This is recorded in the statement. $\Box$

\begin{lemma}\label{2u2}  Let $G^*\in \{ Spin_{2n+1}(q),$  $Spin_{2n}^+(q),Spin_{2n}^-(q), Sp_{2n}(q)\}$, q odd. Let $a,b,c\geq 0$ be integers
such that $a+b+c=n$, $a\neq 1,b\neq 1$,  and if $G^*=Spin_{2n}^-(q)$ then $b+c<n$.
If $G^*\neq Sp_{2n}(q)$ then suppose that $b(q-1)/2$ is even.
Suppose that $q\geq n-a-b+5$. Then there exists a semisimple element $s\in G^*$ such that

\medskip
\begin{center}
$\nu(C_{G^*}(s))=\begin{cases}
\al(a)\al(b)\beta(c) &\text{if } ~G^*=Sp_{2n}(q);\\
\al(a)\al^+(b)\beta(c) &\text{if } ~G^*=Spin_{2n+1}(q);\\
\al^+(a)\al^+(b)\beta(c) &\text{if } ~G^*=Spin^+_{2n}(q);\\
\al^-(a)\al^+(b)\beta(c)&\text{if } ~G^*=Spin^-_{2n}(q).
\end{cases}$
\end{center}
\end{lemma}

Proof. Let $V$ be the natural module for $G^*$ and the respective classical group.  Consider an orthogonal decomposition $V=W_1\oplus W_2\oplus W_3$, where $\dim W_1=2a$ or $2a+1$,
$\dim W_2=2b$ and $W_3=(W_1+W_2)^\perp$. If $V$ is orthogonal, choose both $W_2,W_3$ to be of Witt defect 0. The condition $b+c<n$  makes this possible if $G^*=Spin_{2n}^-(q)$, in the other cases
this is well known to be possible.

Choose a basis $B=\{b_0,b_1\ld b_{2n}\}$ in $V$, where $b_0$ is dropped unless $G^*=Spin_{2n+1}(q)$.
We can assume that $b_{2n-2c+1}\ld b_{2n}\in W_3$, $b_{2n-2c-2b+1}\ld b_{2n-2c}\in W_2$ and the remaining elements of $B$ are in $W_1$. With respect to this basis consider the matrix  $t=\diag(\Id,-\Id_{2b},s')$,
where $s'$ is in $\Omega_{2c}^+(q)$ or $Sp_{2c}(q)$. (Note that $-\Id_{2b}\in\Omega(W_2)$ by Lemma \ref{kk2}.) By Lemma~\ref{rr1}, we can choose $s'$ to be such that $\pm 1$ are not \eis of $s'$ 
and the number of unipotent characters of $C_{SO(W_3)}(s')$ or
$C_{Sp(W_3)}(s')$  equals
$\be(c)$. 
If $G^*=Sp_{2n}(q)$ then $C_{G^*}(t)=Sp(W_1)\times Sp(W_2)\times C_{Sp(W_3)}(s')$, and the result follows   as $\nu(Sp(W_1))=\al(a)$ and $\nu(Sp(W_2))=\al(b)$.

Suppose that $G^*$ is orthogonal. Then $\nu(SO(W_2))=\al^+(b)$ as $W_2$ is of Witt defect 0,
whereas $\nu(SO(W_1))=\al(a)$, $\al^+(a)$ or $\al^-(a)$ depending on whether $G^*=Spin_{2n+1}(q),$  $Spin_{2n}^+(q),Spin_{2n}^-(q)$, respectively. So again the result follows from  Lemma~\ref{st5}.  $\Box$


 \section{Some relations between $\al(n)$, $\al^+(n)$, $\al^-(n)$ and $\beta(n)$}

For $x\in{\mathbb R}$ let $[x]$ denote the maximum integer that does not exceed $x$.

The enumeration of unipotent characters in our context has a nice combinatorial description
(see \cite[Theorem 8.2]{Lu} or \cite[Section 13.8]{Ca}); for computing  $\al(n)$, $\al^+(n)$
and  $\al^+(n)$ for small $n$ (see Table 2) we use Lusztig's formulae   \cite[\S 3]{Lu} expressing
these functions in terms of $p(m)$ with $m\leq n$.

\begin{lemma}\label{lu1}
For $n\in \NN$ odd, $\alpha^- (n)= \alpha^+ (n)$, and
for all  $n\in \NN$ $\alpha^- (n)\leq \alpha^+ (n)
\leq \alpha (n)$.
\end{lemma}

Proof. From Lusztig's generating function \cite[(3.4.2)]{Lu}, we have
$\al^+(n)-\al^-(n)=0$ for $n$ odd, and   $\al^+(n)-\al^-(n)=2 p(n/2)$ for $n$ even,
so always $\al^+(n)\ge \al^-(n)$.

Let $p_2(n)$ denote the number of pairs of partitions that sum up to~$n$, hence $p_2(n)=\sum_{m = 0}^n p(m)p(n-m)$. Again from \cite{Lu}, we have
$$
\al^+(n) = \frac 12 p_2(n) + \frac 34 (1+(-1)^n) p(\frac n2) + \sum_{j>1 \text{ even}} p_2(n-j^2)
$$
and
$$\al(n)=\sum_{j\ge 0} p_2(n-(j^2+j))\:.$$
For $n\le 13$ the claim is easily checked directly (see Table 2). For $n\ge 14$, the
easy inequality $\frac 32 p(\frac n2 ) < p(n-6)+p(n-7)< p_2(n-6)$
and a comparison of the summands in the sums above gives the claim. $\Box$


\begin{propo}\label{albe}
For $n\le 43$ we have $\al(n) > \beta(n)$.
For all $n> 43$, we have $\al(n)<\beta(n)$.\end{propo}

Proof.  For $n\le 43$ the stated inequalities for $\al(n)$  hold by computation (see Table 2); these  also show that  $\al(n)<\beta(n)$ for $44\le n \le 300$.

For $n>43$, we use very rough estimates to give an upper bound for $\al(n)$.
First, we have $p(n) < 2^{[n/2] +1}$ for all $n$  (for example, use \cite{BO}).
Hence, $p_2(m) = \sum_{i=0}^m p(i)p(m-i) \leq (m+1) 2^{[m/2] +2}$,
for all $m$. Applying this, we have for any $n\ge 2$:
$$\al(n)=\sum_{j\ge 0} p_2(n-(j^2+j))\le (n^2-1) 2^{[n/2] +2}\:.$$
One easily checks that  $\al(n)\leq(n^2-1) 2^{[n/2] +2} < 5^{(n-3)/4}$ for  $n\ge 244$.
Using Theorem \ref{bo2}
we conclude  that $\al(n)<\be(n)$ for all $n>43$. $\Box$


\begin{propo}\label{albe1}
For $2<n\le 38$, we have $\beta(n)<\al^-(n)\leq \al^+(n) $.
For all $n \ge 39$, we have $\al^-(n)\leq\al^+(n)<\beta(n) $.\end{propo}

Proof.  For  $2<n\le 43$ the stated inequalities hold by Table 2, for $n>43$
these follow from Lemma~\ref{lu1} and Proposition~\ref{albe}. $\Box$


\begin{corol}\label{c21} Let $n>43$. Then for fixed $n$ but varying $a$, the maximum of each function $\al(a)\beta(n-a)$, $\al^+(a)\beta(n-a)$ and $\al^-(a)\beta(n-a)$
is attained for $a\leq 43$.
\end{corol}

Proof. Suppose on the contrary that the maximum is attained at some $a>43$. Then
$\al^-(a)\beta(n-a)\leq \al^+(a)\beta(n-a)\leq \al(a)\beta(n-a)<\beta(a)\beta(n-a)\leq \beta(n)$.
As $n>43$, by Theorem~\ref{bo2}, we have $\beta(n)=5\beta(n-4)<10\beta(n-4)= \al^-(4)\beta(n-4)<\al^+(4)\beta(n-4)<\al(4)\beta(n-4)$, a contradiction.    $\Box$

\smallskip

\subsection{The products $\al(a)\beta(n-a)$, $\al^+(a)\beta(n-a)$ and $\al^-(a)\beta(n-a)$}

\begin{lemma}\label{dtt}
$(1)$ If $13<a\le 43$ then $\al(a+4)<5\al(a)$,
$\al^+(a+4)<5\al^+(a)$ and $\al^-(a+4)<5\al^-(a)$.

$(1a)$ If $0<a\leq 13$, then $\al(a+4)>5\al(a)$, $\al^+(a+4)>5\al^+(a)$ and $\al^-(a+4)>5\al^-(a)$.

$(2)$  If $13<a\le 43$ then
 $\al(a+4)\beta(m)<\al(a)\beta(m+4)$, $\al^+(a+4)\beta(m)<\al(a)^+\beta(m+4)$ and $\al^-(a+4)\beta(m)<\al^-(a)\beta(m+4)$  for every integer $m\geq 0$.

$(3)$  If $a\leq 13$,
$m>3,m\neq 5,6,11$ then $\al(a)\beta(m)<\al(a+4)\beta(m-4)$, $\al^+(a)\beta(m)<\al^+(a+4)\beta(m-4)$ and $\al^-(a)\beta(m)<\al^-(a+4)\beta(m-4)$.

  More precisely, $\al(a)\beta(m)<\al(a+4)\beta(m-4)$ if $m=11,a<13$ or  $m=6,a<12$ or $m=5,a<8$.

$(4)$  If $a<n$ then the maximum of $\al(a)\beta(n-a)$, $\al^+(a)\beta(n-a)$ and $\al^-(a)\beta(n-a)$  is attained for $a\leq 17$. If $n>24$ then, additionally, $a>13$.\end{lemma}

Proof. (1), (1a) follows directly by Table 2.

(2) By Theorem~\ref{bo2} and Table 1, we have $ 5\beta(m)\leq \beta(m+4)$
so the  claim follows from (1). (Note that $5\beta(m)=\beta(m+4)$ for $m\neq 1,2,7$.)

(3) If $m>3,m\neq 5,6,11$ then $5\beta(m-4)=\beta(m)$.
Therefore, $\al(a+4)\beta(m)>5\al(a)\beta(m)=\al(a)\beta(m-4)$ by (1a). Similarly for $\al^+(a), \al^-(a)$ in place of $\al(a)$.

 Let $m=11$. Then $\beta(11)=77$, $\beta(7)=15$, so $\beta(11)=\frac{77}{15}\beta(7)$. So the result follows if   $\al(a+4)>\frac{77}{15}\al(a)$.
  This is true if $a<13$.

Let $m=6$. Then $\beta(6)=11$, $\beta(2)=2$, so $\beta(6)=\frac{11}{2}\beta(2)$. So
the result follows if   $\al(a+4)>\frac{11}{2}\al(a)$. This is true for $a<12$.

Let $m=5$. Then $\beta(5)=7$, $\beta(1)=1$, so $\beta(5)=7\beta(1) $. So
the result follows if   $\al(a+4)>7\al(a)$. This is true for $a<8$.

(4) By Corollary~\ref{c21} we may assume that $a\le 43$.
Suppose that $a>17$. Then, by (2), $\al(a)\beta(n-a)<\al(a-4)\beta(n-a+4)$, a contradiction.
If $a\leq 13$   then $n>24$  implies $n-a> 11$, so  $\al(a)\beta(n-a)<\al(a+4)\beta(n-a-4)$ by (3). Similarly for $\al^+(a), \al^-(a)$ in place of $\al(a)$.  $\Box$


 \begin{propo}\label{q-even}
$(1)$ For $n < 18$ the maximum of
$ \al(a)\beta(n-a)$,  $\al^+(a)\beta(n-a)$ and $\al^-(a)\beta(n-a)$
is attained for $a=n$.

$(2)$ Let $n\geq 18$. Then the maximum of $\al(a)\beta(n-a)$
is attained for $a=16,15,14,15$ when $n\equiv 0,1, 2,3 \mod 4$.

$(3)$ The maximum of
$\al^+(a)\beta(n-a)$ and of $\al^-(a)\beta(n-a)$
is attained for $a= 16,17,14,15$ when $n\equiv 0,1,2,3 \mod 4$, respectively (in particular, $n\equiv a\mod 4$).\end{propo}

Proof. By computer calculation the claim is easily checked up to $n=29$.
Let $n>29$. By Lemma~\ref{dtt}(4), the maximum of each of these functions is attained for $a$ with $13<a\leq 17$. Then $n-a>7$. Write $n-a=r+4k$, where $7<r\leq 11$ and $k\ge 0$ is an integer.
By Theorem~\ref{bo2}, $\beta(n-a)=5^k\beta(r)$. So
$\al(a)\beta(n-a)=5^k\al(a)\beta(r)$, where $a+r<29$. By the above, the maximum of
$\al(a)\beta(r)$ is attained for $a= 16,15,14,15$  if $a+r$ is congruent to $0,1,2,3 \mod 4$, respectively. Say, if $4|(a+r)$ then $a=16$, and $\al(a)\beta(n)=5^k\al(a)\beta(r)\leq 5^k\al(16)\beta(r)= \al(16)\beta(n)$, whence the result. The other cases are similar, as well as the cases with $\al^+(a)$, $\al^-(a)$ in place of $\al(a)$. $\Box$


\subsection{The products $\al(a)\al(b)\beta(n-a-b)$, and $\al(a)\al^+(b)\beta(n-a-b)$}

\begin{lemma}\label{ma1} Let
$n,a,b\geq 0$ be integers such that $a+b\leq n$. For $n$ fixed, the maximum of
$\al(a)\al(b)\beta(n-a-b)$, $\al(a)\al ^+(b)\beta(n-a-b)$, $\al^+(a)\al ^+(b)\beta(n-a-b)$
and $\al^+(a)\al ^-(b)\beta(n-a-b)$ is attained for  $a\leq17,b\leq17$.

Furthermore,

$(1)$ if $n-a-b>11$, then $a> 13$, $b> 13$;

$(2)$ if $n>45$,  then $a> 13$, $b> 13$.\end{lemma}

Proof.  The first statement follows from Corollary~\ref{c21} and Lemma~\ref{dtt}(4).
Furthermore,   $n-a>11$ and  $n-b>11$. Suppose that $a\leq 13$. Then, by Lemma~\ref{dtt}(2), we have  $ \al(a)\al (b)\beta(n-a-b)<\al(a+4)\al (b)\beta(n-a-b-4)$, a contradiction. Similarly, for the other three functions, as well as for $b\leq 13$, whence (1).
In addition, if $n>45$ then $n-a-b>11$ as $a+b\leq 34$, whence (2). $\Box$


  \begin{propo}\label{alaalb1}
  Let $n,a,b\geq 0$ be integers, $n\geq a+ b$.

In the table below, we record for each of the functions
$\al(a)\al(b)\beta(n-a-b)$, $\al(a)\al ^+(b)\beta(n-a-b)$, $\al(a)^+\al ^+(b)\beta(n-a-b)$,
and $\al(a)\al ^-(b)\beta(n-a-b)$ the pairs $(a,b)$ where the functions attain their maximum, for $n\ge 28$ or $n\ge 29$; for the first and third function, we list the pairs with $a\ge b$.
\\
(Here we write $\equiv_4$ for the congruence modulo~$4$.)
\medskip

\begin{center}
{
\begin{tabular}{|l|l|l|l|l|l|}\hline
~~~~~~~~~~~~~~~~~~~~~~~~~~\qquad {\rm function} & {\rm bound}
& $n \equiv_4 0 $ & $n \equiv_4 1 $ & $n \equiv_4 2 $& $n \equiv_4 3 $ \\
\hline
$\al(a)\al(b)\beta(n-a-b)$& $n\ge 28$ & $(14,14)$&$(15,14)$&$(15,15)$&$(16,15)$\\
\hline
$\al(a)\al^+(b)\beta(n-a-b)$ &$n\ge 29$ & $(16,16)$&$(15,14)$&$(14,16)$&$(15,16)$\\
\hline
$\al^+(a)\al^+(b)\beta(n-a-b)$& $n\ge 29$ &$ (16,16)$&$(15,14)$&$(16,14)$&$(16,15)$\\
\hline
$\al^+(a)\al^-(b)\beta(n-a-b)$ &$n\ge 29$ & $(16,16)$&$(14,15)$&$(14,16)$&$(16,15)$\\
\hline
\end{tabular}}
\end{center}
\end{propo}
\smallskip

Proof. The assertion was checked to hold for $n \le 50$ by computer, so we may assume $n>50$.
 Lemma~\ref{ma1} shows that the values $a,b$ at which all four  products  in question
 attain their  maximum satisfy $13< a,b < 18$.  Write $n=4k+r$, with $a+b+7\le r\le a+b+10$,
 and some integer $k\ge 0$.
Then $r\leq 44$. Let $\gamma_n(a,b)$ stand for any of the functions above. As $r>7$, by  Theorem~\ref{bo2}  we have $\beta(n-a-b)=\beta(r-a-b)\beta(4)^k$, and hence $\gamma_n(a,b)=\beta(4)^k\gamma_{r}(a,b)$. Since $35\le r\le 44$, the claim holds for $r$, and the result follows. $\Box$

\med
 Remarks. (1) For all $n<32$, the maximum of $\al(a)\al(b)\beta(n-a-b)$
is attained for pairs $(a,b)$ such that $a+b=n$.

(2) For $n\leq 33$ the maxima of the functions $\gamma_n(a,b)$ defined in the proof of Proposition \ref{alaalb1} have been calculated by computer and are shown in Tables 3 and~4 at the end of the paper.

\section{Proof of the main results for $q$ even}

In this section $q$ is even and $G^*\in\{Sp_{2n}(q),Spin^\pm_{2n}(q)\cong \Omega^\pm_{2n}(q)\}$. 
For $q$ large enough, we  determine the maximum of $\nu(C_{G^*}(s))$ when $s$ runs over the semisimple elements of~$G^*$.

Let  $V$ be a vector space of dimension $2n$ over $\FF_q$ viewed as the natural module for $G^*$,
so $V$ is endowed with a suitable form defining $G^*$.
 Denote by $V_1$ the 1-eigenspace of $s$ on $V$. By Lemma~\ref{ei2},  $V_1$ is non-degenerate
and  $\dim V_1=2a$ is even. Set $W=V_1^\perp$, so $V=V_1\oplus W$. Let $s'$ denote the restriction of $s$ to $W$. We keep these notations until the end of this section.

\bl{n18} Proposition $\ref{p11}$ is true for $q$ even.
\el

Proof. Suppose that $G^*\cong Spin^\pm_{2n}(q)$ (the proof for $G^*=Sp_{2n}(q)$ is similar, and hence omitted).
By Lemma \ref{st5}, $\nu(C_{G^*}(s))=\nu(SO(V_1))\cdot \nu(C_{SO(W)}(s'))=\al^\pm (a)\cdot \nu(C_{SO(W)}(s'))$. By Lemma \ref{ei1}(1), $\nu(C_{SO(W)}(s'))\leq \be(n-a)$. If $n<18$ then the maximum of  $\al^\pm (a)\be(n-a)$
is attained for $a=n$ (Proposition \ref{q-even}(1)), which is realized for $s=1$. $\Box$

\subsection{Symplectic groups in even characteristic}

 Let $G^*=Sp_{2n}(q)$. Then $C_{G^*}(s)\cong Sp_{2a}(q)\times C_{Sp_{2(n-a)}(q)}(s')$,
 so $\nu(C_{G^*}(s))=\al(a)  \cdot \nu(C_{Sp_{2(n-a)}(q)}(s'))$, and we are to determine the maximum of this product. Recall that $\nu(C_{Sp_{2(n-a)}(q)}(s'))\leq \be(n-a)$ by Lemma~\ref{ei1}, as $s'$ does not have \ei 1.

\begin{theo}\label{sp2} 
Theorem $\ref{t21}$ is true for $G^*=Sp_{2n}(q)$, q even.
\end{theo}

Proof.  
If $n$ is fixed and $a$ varies, the maximum of $\al(a)\beta(n-a)$ is attained for 
$a=16,15,14,15$ for $n\equiv 0,1,2,3\pmod 4$, respectively,
 see Proposition \ref{q-even}.  So $f(n)=\al(a)\be(n-a)$ for these values of $a$.

 The values of $\al(a)$ for $a\leq 17$
are provided by Table 2, and $\beta(n-a)$ is given by Table 1. Using this, the result follows
by easy computations. (For instance, $n\equiv 1\mod 4$ 
implies $n-15\equiv  2\mod 4$, and then $\be(n-15)=11\cdot 5^{n-21}$ by Table~1.)

Now we turn to the additional statement on the bound being attained for large~$q$.
By Lemma~\ref{2u1}, if $n-a\leq q-5$ then there exists a semisimple  element $s\in G^*$ such that $\nu(C_{G^*}(s))=\al(a)\beta(n-a)$.  This holds if $n-14\leq q-5$, that is, $q\geq n-9$. $\Box$


\subsection{Orthogonal groups in even characteristic}

In this case we are to consider the groups $G^*=Spin^\pm_{2n}(q)\cong SO^\pm_{2n}(q)$. 
By Lemma~\ref{st5}, $\nu(C_{G^*}(s))=\nu(SO(V_1))\cdot \nu(C_{SO(W)}(s'))$, and
  $ \nu(C_{SO(W)}(s'))\leq \beta(n-a)$ by Lemma~\ref{ei1}. So  $|{\mathcal E}_s|=\nu(C_{G^*}(s))\leq \al^\pm(a)\beta(n-a)$, where one chooses $+$ \ii $V_1$ is of Witt defect~0.

\begin{theo}\label{ore2} Theorem $\ref{t21}$ is true if $G^*=Spin^\pm_{2n}(q)$, $q$ even. More precisely,
  $|{\mathcal E}_s|=\nu(C_{G^*}(s))\leq
\al^\pm (a)\beta(n-a)$, where   $a= 16,17,14,15$ when $n\equiv 0,1,2,3 \mod 4$, respectively, and 
$|{\mathcal E}_s|\leq f^\pm(n)$, where $f^\pm(n)$ are as in Theorem $\ref{t21}$.
\end{theo}

Proof. 
Let $x,y$ denote the Witt defect of $V_1$, $W$, respectively.
\smallskip

(i) $(x,y)=(0,0)$. Then  $G^*=Spin^+_{2n}(q)$ and $\nu(C_{G^*}(s))\leq \al^+(a)\beta(n-a)$;
\smallskip

(ii) $(x,y)=(1,1)$. Then $G^*=Spin^+_{2n}(q)$ and $\nu(C_{G^*}(s))\leq \al^-(a)\beta(n-a)$;
\smallskip

(iii) $(x,y)=(1,0)$. Then $G^*=Spin^-_{2n}(q)$ and $\nu(C_{G^*}(s))\leq \al^-(a)\beta(n-a)$;
\smallskip

(iv) $(x,y)=(0,1)$. Then $G^*=Spin^-_{2n}(q)$ and $\nu(C_{G^*}(s))\leq \al^+(a)\beta(n-a)$  if $n-a$ is odd, and $\nu(C_{G^*}(s))\leq\al^+(a)\beta'(n-a)$ if $n-a$ is even, see Lemma~\ref{ei2a}
(and Lemma~\ref{nk8} for the values of $\beta'$).

 Assume  that  $s$ is chosen    so that $\nu(C_{G^*}(s))$ is maximal. Then we show that
cases (ii),(iv) can be ignored for our purpose.

In case (ii), $\al^-(a)\leq \al^+(a)$, and, by Lemma~\ref{2u1},   $\nu(C_{G^*}(s))= \al^+(a)\beta(n-a)$ for some $s\in G^*$. So we do not need to care whether the same maximum can be  attained in case (ii). Then, using Proposition  \ref{q-even}, we obtain the data for $f^+(n)$.
For instance, if $n \equiv 1 \mod 4$ then $a=17$ and $\al^+(17)\beta(n-17)=6007\cdot 5^{(n-17)/4}$.

Suppose, on the contrary, that (iv) holds. We first obtain an upper bound for $\max_{a:\  n-a \text{ even}}\,\al^+(a)\beta'(n-a) $ and $\max_{a:\ n-a \text{ odd}}\,\al^+(a)\beta(n-a) $, and next show that these are less than max$ _a\,\al^-(a)\beta(n-a)$, which will yield the stated claim.

Assume first that  $n-a$ is even. By Lemma~\ref{nk8}, we have
$ \beta'(n-a)=539\cdot \beta(n-a-16)$ if $n-a \equiv 0 \mod 4$, $n-a\geq 16$, and $\beta'(n-a)=49\cdot \beta(n-a-10)$ if $n-a\equiv 2 \mod 4$, $n-a>6$. So $\al^+(a)\beta'(n-a)=539\cdot \al^+(a)\beta(n-a-16)$ and
  $49\cdot \al^+(a)\beta(n-a-10)$ accordingly; by Proposition \ref{q-even}, if $n-16\geq 18$, respectively, $n-10\geq 18$, then the maximum of these is attained when $a=16$  if $n\equiv n-16\equiv 0 \mod 4$, respectively, when $a=14$  if $n-10 \equiv 0 \mod 4$ equivalently, $n \equiv 2 \mod 4$. Thus, if $n-a$ is even, then $\al^+(a)\beta'(n-a)$ does not exceed

\smallskip

{\small
$\begin{cases}
539\al^+(16)\beta(n-32)=539\cdot 4110\cdot 5^{(n-32)/4}&\text{if }  n-a\equiv 0,~n \equiv 0 \mod 4;\\
539\al^+(17)\beta(n-33)=539\cdot 6007\cdot 5^{(n-33)/4}&\text{if }  n-a\equiv 0,~n \equiv 1 \mod 4;\\
539\al^+(14)\beta(n-30)=539\cdot 1836\cdot 5^{(n-30)/4}&\text{if }  n-a\equiv 0,~n \equiv 2 \mod 4;\\
539\al^+(15)\beta(n-31)=539\cdot 2730\cdot 5^{(n-31)/4}&\text{if }  n-a\equiv 0,~n \equiv 3 \mod 4;\\
49\al^+(16)\beta(n-26)=539\cdot 4110\cdot 5^{(n-32)/4}&\text{if }   n-a\equiv 2,~n \equiv 0 \mod 4;\\
49\al^+ (17)\beta(n-27)=539\cdot 6007 \cdot 5^{(n-33)/4}&\text{if } n-a\equiv 2,~n \equiv 1 \mod 4;\\
49\al^+ (14)\beta(n-24)=539\cdot 1836  \cdot 5^{(n-30)/4}&\text{if }  n-a\equiv 2,~n \equiv 2 \mod 4;\\
49\al^+ (15)\beta(n-25)=539\cdot 2730  \cdot 5^{(n-31)/4}&\text{if }  n-a\equiv 2,~n \equiv 3 \mod 4.
\end{cases}$}

\smallskip
\noindent
Indeed, here $n-i\equiv 0 \mod 4$ for $i=32,33,30,31$ in the last four rows, so $\beta(n-i+6)=\be(6)\beta(n-i)=11\be(n-i)$ for these $i$ by Theorem~\ref{bo2}, whence the equalities there as $49\cdot 11=539$.
(Note that we do not assert that these bounds are attained.)

 Next we assume  $n-a$ to be odd in  case (iv). Then $\beta(n-a)=7\beta(n-5-a)$ if $n-a\equiv 1 \mod 4$ and $\beta(n-a)=77\beta(n-11- a) $ if $n-a\equiv 3 \mod 4$, see Theorem~\ref{bo2}.
By Proposition  \ref{q-even},  applied to  $\al^+(a)\beta(n-a-5)$ and  $\al^+(a)\beta(n-a-11)$,  if $n-11\geq 18$ then  $\al^+(a)\beta(n-a)$ with $n-a$ odd does not exceed the following values:

\smallskip

{\small
$\begin{cases}
7\al^+(16)\beta(n-21)=7 \cdot 4110\cdot 5^{(n-21)/4}&\text{if }  n-a\equiv 1,~n \equiv 1 \mod 4;\\
7\al^+(17)\beta(n-22)=7 \cdot 6007\cdot 5^{(n-22)/4}&\text{if }  n-a\equiv 1,~n \equiv 2 \mod 4;\\
7\al^+(14)\beta(n-19)=7 \cdot 1836\cdot 5^{(n-19)/4}&\text{if }  n-a\equiv 1,~n \equiv 3 \mod 4;\\
7\al^+(15)\beta(n-20)=7 \cdot 2730\cdot 5^{(n-20)/4}&\text{if }  n-a\equiv 1,~n \equiv 0 \mod 4;\\
77\al^+(16)\beta(n-27)=77 \cdot 4110\cdot 5^{(n-27)/4}&\text{if } n-a\equiv 3 ,~n \equiv 3 \mod 4;\\
77\al^+ (17)\beta(n-28)=77 \cdot 6007\cdot 5^{(n-28)/4}&\text{if }  n-a\equiv 3,~n \equiv 0 \mod 4;\\
77\al^+ (14)\beta(n-25)=77 \cdot 1836\cdot 5^{(n-25)/4}&\text{if }  n-a\equiv 3,~n \equiv 1 \mod 4;\\
77\al^+ (15)\beta(n-26)=77 \cdot 2730\cdot 5^{(n-26)/4}&\text{if }  n-a\equiv 3,~n \equiv 2 \mod 4.
\end{cases}$}
\smallskip

We have to compare this with $\max_a \al^-(a)\beta(n-a)$. By Proposition  \ref{q-even}, this is equal to

\smallskip
$\begin{cases}
\al^-(17)\beta(n-17)=6007\cdot 5^{(n-17)/4}=  6007\cdot 25\cdot 5^{(n-25)/4}&\text{if }  n \equiv 1 \mod 4;\\
\al^-(15)\beta(n-15)=2730\cdot 5^{(n-15)/4}=  2730\cdot 5^{3}\cdot 5^{(n-27)/4}&\text{if }  n \equiv 3 \mod 4;\\
\al^-(16)\beta(n-16)=4066\cdot 5^{(n-16)/4}=  4066\cdot 5^{3}\cdot 5^{(n-28)/4}&\text{if }  n \equiv 0 \mod 4;\\
\al^-(14)\beta(n-14)=1806\cdot 5^{(n-14)/4}=  1806\cdot 5^{3}\cdot 5^{(n-26)/4}&\text{if }  n \equiv 2 \mod 4.
\end{cases}$

\smallskip\noindent
Then we  conclude that the latter are greater than the former.
So in case (iv) with $n\geq 34$ the maximum of $\nu(C_{G^*}(s))$ does not exceed
$\al^-(a)\beta(n-a)$.  The same trivially holds in case (iii). So the values
for $f^\pm (n)$ for $n\geq 34$ follow from the above.

For $n<34$ we use computer calculations.

Finally, we show that the bound $|{\mathcal E}_s|\leq f^\pm (n)$ is attained for $q$ as stated.
We have $13<a\leq 17$, so $n-a\leq n-14$ above. By Lemma~\ref{2u1},
if $n-a\leq q-5$ then there exists a semisimple  element $s\in G^*$ such that $\nu(C_{G^*}(s))=\al^\pm (a)\beta(n-a)$.
This holds if $n-14\leq q-5$, that is, $q\geq n-9$. So the result follows. $\Box$

\section{Proof of the main results for $q$ odd}

In this section we assume that  $q$ is odd. Let $V$ be the natural $\FF_qG^*$-module, and
$s\in G^*$  a semisimple element. Let $V_1$ and $V_2$ denote the 1- and $-1$-eigenspaces of $s$ on $V$, respectively. These spaces are non-degenerate (if non-zero), and have even dimensions, except for the case where $\dim V$ is odd  (see Lemma~\ref{ei2}).
Set   $\dim V_1=2a$ or $2a+1$ and  $\dim V_2=2b$, where $0\leq a,b\leq n$. Set $W=(V_1 +V_2)^\perp$. Then $V=V_1\oplus V_2\oplus W$. One easily observes that $C_{G^*}(s)$ stabilizes $V_1, V_2$ and $W$. Let $s'$ be the restriction of $s$ to $W$. As above,  $|\mathcal{E}_s|=\nu(C_{G^*}(s))$.

\bl{pp2} Proposition $\ref{p11}$ is true for $q$ odd.
\el

Proof. Tables 3, 4 at the end of the paper are obtained by computer calculations,
and give us the maximum of the functions in question. So we have to show that these coincide with the maximum of  $|{\mathcal E}_s|$ in each case.

From a look at the tables, one observes that  $n\leq 32$ implies $n=a+b$ (in the notation of the tables). Let $V$ be the natural module for $G^*.$ Write $V=V_1+V_2$, where $V_1,V_2$ are non-degenerate subspaces of $V$ such that $\dim V_1=2a$ or $2a+1$, $\dim V_2=2b$ and $V_2$ is of Witt defect 0. This is always possible  unless $G^*=Spin^-_{2n}(q)$ and $a=0 $; this happens only for $n=2,4,6$
 (see Table 4),  but these cases are excluded in the  statement.
Consider the matrix $t=\diag(\Id_{2a},-\Id_{2b})\in SO(V)$ such that $V_1,V_2$ are
the 1- and $-1$-eigenspace of $t$ on $V$. If $b=0$ then we set $t=\Id_{2a}$ or $\Id_{2a+1}$.
From another look at the tables we conclude that $b$ is always even, and hence  $-\Id\in \Omega^+_{2b}(q)$ for $b>0$ by Lemma~\ref{kk2}. Then we take $s$ from the preimage of $t$ in $G^*$. $\Box$

\med
Remark. In case of $SO_{2n}^-(q)$ and $n=2,4,6$, one easily checks that the maximum of $|{\mathcal E}_s|$ is attained for $\al^-(2)\al^+(0)=2$,  $\al^-(4)\al^+(0)=10$ and $\al^-(4)\al^+(2)=40$ respectively. If $G^*=Sp_{64}(q)$ then $\max_{s\in G^*}|{\mathcal E}_s|=5\al(14)^2$ by Theorem~\ref{sp88}.

\subsection{$G^*$ is symplectic}

\begin{theo}\label{sp88}  Theorem $\ref{t21}$ is true if 
 $ G^*\cong Sp_{2n}(q)$, $q$ odd.  
More precisely, $|{\mathcal E}_s|=\nu(C_{G^*}(s)) \leq \tau(n)=\max_{a+b\le n}\, \al(a)\al(b)\beta(n-a-b)$, where
$\tau(n)$  is as in Theorem $\ref{t21}$.
\end{theo}

  Proof. One easily observes that $C_{G^*}(s)$ stabilizes $V_1, V_2$ and $W$. 
  So    $C_{G^*}(s)\subset Sp_{2a}(q)\times Sp_{2b}(q)\times Sp_{2(n-a-b)}(q) $, and in fact $C_{G^*}(s)= Sp_{2a}(q)\times Sp_{2b}(q)\times C_{Sp_{2(n-a-b)}(q)}(s') $. Therefore, $\nu(C_{G^*}(s))=\al(a)\al(b)
\cdot\nu(C_{Sp(W)}(s'))$. By Lemma~\ref{ei1}, $\nu(C_{Sp(W)}(s'))\leq \beta(n-a-b)$. 
 Therefore, if $s$ varies,  $|\mathcal{E}_s|$
does not exceed the maximum of $\al(a)\al(b)\beta(n-a-b)$, where $a,b\geq 0$ and  $a+b\leq n$.
The values of $a,b$
for which the function $\al(a)\alpha (b)\beta(n-a-b) $ attains its maximum
are determined in Table 3 for $n\le 33$ and in Proposition \ref{alaalb1} for $n\geq 28$.
 This yields the explicit expressions for $\tau(n)$ 
in Theorem $\ref{t21}$.

It remains to show the additional statement in Theorem \ref{t21}. By Lemma \ref{2u2}
(or Lemma~\ref{rr1} applied to $Sp_{2(n-a-b)}(q)$),  if $q\geq n-a-b+5$ with $a,b$  as above, there is a semisimple  element $s\in G^*$ such that
$|\mathcal{E}_s|=\nu(C_{G^*}(s))=\al(a)\alpha (b)\be(n-a-b)$. So the bound is attained, whence the result.  $\Box$

\smallskip

\subsection{The case $G^*=Spin_{2n+1}(q)$  }

\begin{theo}\label{od1} Theorem $\ref{t21}$ is true for
  $ G^*=Spin_{2n+1}(q)$, $q$ odd.
 More precisely,  $|{\mathcal E}_s|=\nu(C_{G^*}(s)) \leq \theta(n)=\max_{a+b\le n}\, \al(a)\al(b)\beta(n-a-b)$, where
$\theta(n)$  is as in Theorem $\ref{t21}(4)$.
\end{theo}

  Proof.  Suppose for a moment that $s$ is an arbitrary semisimple element in $G^*$.
Then $|{\mathcal E}_s|=\nu(SO( V_1))\cdot \nu(SO( V_2))\cdot \nu(C_{SO(W)}(s'))=\al(a)\al^\pm (b)\nu(C_{SO(W)}(s'))$ by Lemma~\ref{st5}, where one chooses the sign $+$ \ii the Witt defect of $W$ is~1. Therefore, $|{\mathcal E}_s|\leq \max_{a,b}\al(a)\al^+ (b)\be(n-a-b)$.
   By  Proposition   \ref{alaalb1}, if $n\geq 32$ then the maximum of $\al(a)\al^+(b)\beta(n-a-b)$ is attained for $(a,b)$   given in the table there (in particular, with $b$ even).  The data $\al(i)\al^+(j)$ for $14\leq i,j\leq 16$ follows from Table~3. So the inequality $|{\mathcal E}_s| \leq \theta(n)$ follows.

For the additional statement in Theorem \ref{t21} we note that the existence of $s$ such that $|{\mathcal E}_s| = \theta(n)$ follows from Lemma~\ref{2u2} (provided $b$ is even, which is the case here).   $\Box$

\smallskip

\subsection{Orthogonal groups of even dimension}

In this subsection we assume that $q$ is odd and $G^*=Spin^\pm_{2n}(q)$.

By Lemma~\ref{st5},     $|{\mathcal E}_s|=\nu(C_{G^*}(s))=\nu(SO( V_1))\cdot \nu(SO( V_2))\cdot \nu(C_{SO(W)}(s'))$. As  $\al^-(a)\leq \al^+(a)$ and $\nu(C_{SO(W)}(s'))\leq \beta(n-a-b)$ (Lemma~\ref{ei1}(1)), it follows that
$|{\mathcal E}_s|\leq max_{a,b}\, \al^+(a)\al^+(b)\be(n-a-b)$. In turn, by Proposition \ref{alaalb1}, if $n \ge 32$ then the maximum of $\al^+(a)\al^+(b)\be(n-a-b)$ is attained for $(a,b)=(16,16),(15,14),(16,14),(15,16)$  for $n\equiv 0,1,2,3\pmod 4$, respectively. In particular,  $b\in\{14,16\}$ is even and $n-a-b\equiv 0\pmod 4$.

\begin{theo}\label{ev+} 
  Theorem $\ref{t21}$ is true for $G^*=Spin_{2n}^+ (q)$, $q$ odd.
More precisely,  $|{\mathcal E}_s|=\nu(C_{G^*}(s)) \leq \theta^+(n)=\max_{a+b\le n}\, \al^+(a)\al^+(b)\beta(n-a-b)$, where
$\theta^+(n)$  is as in Theorem $\ref{t21}(5)$.
\end{theo}

Proof.  The comments prior to the theorem show that $|{\mathcal E}_s| \leq \theta^+ (n)$, so we are left to show that the equality holds for some semisimple  element $s\in G^*$. This follows from Lemma~\ref{2u2} as now $b\in\{14,16\}$ and $b(q-1)/2$ are even, so Lemma~\ref{2u2} applies. $\Box$


\begin{theo}\label{ev-} Theorem $\ref{t21}$ is true for
  $s\in G^*=Spin_{2n}^- (q)$, $q$ odd. 
\end{theo}

Proof.  Let $V$  be the natural module for $G^*$ and $s\in G^*$ an arbitrary semisimple element.   As above, consider a decomposition
$V=V_1\oplus V_2\oplus W$, where  $V_1$ and $V_2$ are the 1- and $-1$-eigenspaces of $s$ on $V$ and $W=( V_1+V_2)^\perp$. By Lemma~\ref{st5}, $|{\mathcal E}_s|=\nu(C_{G^*}(s))=\nu(SO(V_1))\cdot \nu(SO(V_2))\cdot \nu(C_{SO(W)}(s'))$, where $s'$ is the restriction of $s$ to $W$.
As the Witt defect of $V^*$ equals 1, for the Witt defects of  $V_1, V_2$ and $W$ we have the following options:

(i) the Witt defect of $V_1$ is 1, the two other are 0;

(ii)  the Witt defect of $V_2$ is 1, the two other are 0;

(iii)  the Witt defect of $W$ is 1, the two other are 0;

(iv)   $V_1,V_2,W$ are of  Witt defect  1.

Suppose that $s$ is chosen so that $|{\mathcal E}_s| $ is maximal. Then  (iv)
can be ignored. Indeed, in this case  $\nu(C_{G^*}(s))=\al^-(a)\al^-(b)\cdot \nu(C_{SO(W)}(s'))$, where $\nu(C_{SO(W)}(s'))\leq \beta(n-a-b)$ by Lemma~\ref{ei1}.  One can choose another  element $s_1\in G^*$ for which
the $-1$-eigenspace is the same as for $s$, the 1-eigenspace $U$, say, is of dimension $2a$ and of Witt defect 0,
and for $W_1=(U+ V_2)^\perp$ choose $s_1'\in \Omega(W_1)$ so that $\nu(C_{SO(W)}(s_1'))= \beta(n-a-b)$. This is possible as the Witt defect of $W_1$ is 0, so $\Omega(W_1)\cong \Omega_{2(n-a-b)}^+(q)$, see Lemma~\ref{rr1}. Then  $|{\mathcal E}_{s_1}| =\nu(C_{G^*}(s_1))=\nu(SO(V_1))\cdot \nu(SO(V_2))\cdot \nu(C_{SO(W)}(s_1'))=\al^+(a) \al^-(b)\cdot \nu(C_{SO(W)}(s_1'))=\al^+(a)\al^-(b)\beta(n-a-b)$.
As  $\al^-(a)\leq \al^+(a)$, we have $|{\mathcal E}_s|\leq |{\mathcal E}_{s_1}| $, so we can assume that (iv) does not hold.

Suppose first that  $W$ is of Witt defect 0.  Then (i) or (ii) holds, and $\nu(C_{G^*}(s))= \al^-(a) \al^+(b)\nu(C_{SO(W)}(s'))$ in case (i) and
$\nu(C_{G^*}(s))= \al^+(a) \al^-(b)\nu(C_{SO(W)}(s'))$ in case (ii). By Lemma~\ref{ei1}(1), $\nu(C_{SO(W)}(s'))\leq \be(n-a-b)$.  By Proposition~\ref{alaalb1},
the maximum of the function $\al^-(a) \al^+(b)\be(n-a-b)$ is attained for $(a,b)=(16,16),(15,14),(16,14),(15,16)$, for $n\equiv 0,1,2,3\pmod 4$, respectively,  and the maximum of  $\al^+(a)\al^-(b)\be(n-a-b)$ is attained for $(a,b)=(16,16),(14,15),(14,16)$, $(16,15)$ for $n\equiv 0,1,2,3\pmod 4$. By Lemmas \ref{ei1}(3)  and \ref{rr1}, if $q>n-a-b+5$ then there is $s'$ such that $\nu(C_{SO(W)}(s'))=\be(n-a-b)$; for the above values of $a,b$ it suffices to assume $q>n-27$.

In case (i),
$b\in\{14,16\}$ is even and the Witt defect of $V_2$ is 0; so $-\Id_{2b}\in\Omega(V_2)$ by Lemma~\ref{kk2}. \itf $t\in\Omega_{2n}^-(q)$, where  $t=\diag(\Id_{2a},-\Id_{2b}, s')$ and $s'\in \Omega(W)$ is such
that $\pm 1$ are not \eis of $s$ and $\nu(C_{SO(W)}(s_1'))=\be(n-a-b)$.
Let $s\in G^*$ be such that $t$ is the matrix of $s$ on $V$. Then $\nu(C_{G^*}(s))=\al^-(a)\al^+(b)\cdot\be(n-a-b)$, where $(a,b)$ are as above. Therefore, in case (i) the maximum of $\nu(C_{G^*}(s))$ is attained for the values of $a,b$ as in the statement, and hence it  is left to be shown that the maximum is not greater than this in cases (ii),(iii).

Suppose that $n$ is odd. Then the maximum of $\al^+(a) \al^-(b)\be(n-a-b)$ in case~(ii) and of $\al^+(a) \al^+(b)\be(n-a-b)$ in case (iii) is attained for $(a',15)$ or $(15,b')$ with $a',b'$ even; in addition, $ \al^-(15)= \al^+(15)$. By swapping $V_1,V_2$ if necessary, we arrive at the case with $a=15$ and $V_1$ of Witt defect 1, that is, at case~(i).  So
the result follows for $n$ odd.

Let $n$ be even. We show that case (ii) can be ignored. Indeed,  the maximum of the functions
 $\al^+(a) \al^-(b)\be(n-a-b)$ and  $\al^-(a) \al^+(b)\be(n-a-b)$ is attained for  $(a,b)=(16,16)$ if $n\equiv 0 \mod 4$, and for $(a,b)=(16,14)$ and $ (14,16)$,
if $n \equiv 2 \mod 4$. So the two maxima coincide.

  It remains to compare the maxima of $|{\mathcal E}_s|$ in cases (i) and (iii) for $n$ even. In case  (i), this is  $\al^-(16) \al^+(16)\be(n-32)$ if $n\equiv 0 \mod 4$, and $ \al^-(16) \al^+(14)\cdot\be(n-30)$ if $n \equiv 2 \mod 4$.
In case (iii), $\nu(C_{G^*}(s))= \al^+(a) \al^+(b)\nu(C_{SO(W)}(s'))\leq \al^+(a) \al^+(b)\be(n-a-b)$.

Suppose first that $n-a-b$ is odd. Then $\beta(n-a-b)=7\beta(n-5-a-b)$ if $n-a-b\equiv 1\mod 4$ and $\beta'(n-a-b)=77\beta(n-11- a-b) $ if $n-a-b\equiv 3 \mod 4$, see Theorem~\ref{bo2}.
In the latter case, Proposition  \ref{alaalb1}  applied to
$\al^+(a)\al^+(b)\beta(n-11-a-b)$,  if $n-11\ge 32$,
yields that  $\al^+(a)\al^+(b)\beta(n-11-a-b)$
does not exceed the following values:

\smallskip

$\begin{cases}
77\al^+ (15)\al^+(14)\beta(n-40)=385945560\cdot 5^{(n-40)/4}& \text{if }
~n \equiv 0 \mod 4;\\
77\al^+ (15)\al^+(16)\beta(n-42)=863963100\cdot 5^{(n-42)/4}& \text{if }
~n \equiv 2 \mod 4.
\end{cases}$
\medskip

A similar statement can be written for $\al^+(a)\al^+(b)\beta(n-5-a-b)$, but one observes from Theorem~\ref{bo2} that $\beta(n-5-a-b)=\be(6)\be(n-11-a-b)$ provided $n-11-a-b\geq  0$ and $n-11-a-b\equiv 0\pmod 4$. As $\be(6)=11$, we obtain the same values as above for the  maximum of $\al^+(a)\al^+(b)\beta(n-5-a-b)$.

Note that  $\beta(n-i)= 5^{(n-i)/4}$ for $i=40,42$ above   as $ n-i\equiv 0 \mod 4$ in each case. (Observe that we do not apply Proposition  \ref{alaalb1} directly to the function  $\al^+(a)\al^+(b)\beta(n-a-b)$ as $a+b$ is here odd whereas  the maximum of this function with $n$ even is attained with $a+b$ even.)

Let $n-a-b$  be even. Then $\nu(C_{SO(W)}(s'))\leq \be'(n-a-b)$.  Recall (Lemma~\ref{nk8})  that
$ \beta'(n-a-b)=539\cdot \beta(n-a-b-16)$ if $n-a-b\equiv 0 \mod 4$, $n-a-b\geq 16$, and $49\cdot \beta(n-a-b-10)$ if $n-a-b\equiv 2 \mod 4$, $n-a-b>6$. So $\al^+(a) \al^+(b)\beta'(n-a-b)=539\cdot \al^+(a)\al^+(b)\beta(n-a-b-16)$ and
  $49\cdot \al^+(a)\al^+(b)\beta(n-a-b-10)$ accordingly; by Proposition \ref{q-even}, if $n-16\geq 18$, respectively, $n-10\geq 18$ then the maximum of these functions is attained for $(a,b)=(16,16)$  if $n\equiv n-16\equiv 0 \mod 4$, respectively, for $(a,b)=(16,14)$  if $n-10 \equiv 0 \mod 4$ (i.e., $n \equiv 2 \mod 4$). In fact, it suffices to record an upper bound
for the case where $n-a\equiv 0 \mod 4$ as $\be(n-a-b-10)=\beta(6)\beta(n-a-b-16)=11\be(n-a-b-16)$ provided $n-a-b-16\geq 0$ and $n-a-b-16\equiv 0\pmod 4$.

Thus, if $n-a-b$ is even then $\al^+(a)\al^+(b)\beta'(n-a-b)$ does not exceed the following values:

\smallskip

$\begin{cases}
539\al^+(16)\al^+(16)\beta(n-48)=9104841900\cdot 5^{(n-48)/4}& \text{if }
~n \equiv 0 \mod 4;\\
539\al^+(16)\al^+(14)\beta(n-46)=4067272440\cdot 5^{(n-46)/4}&\text{if }
~n \equiv 2 \mod 4.\\
\end{cases}$

\smallskip
\noindent
These must be compared with  $\al^-(16) \al^+(16)\be(n-32)=16711260\cdot 5^{(n-32)/4}$ if $n\equiv 0 \mod 4$, and $ \al^-(16) \al^+(14)\be(n-30)=7465176  \cdot 5^{(n-30)/4}$ if $n \equiv 2 \mod 4$.


For $n\equiv 0 \mod 4$ we have
$\al^-(16)\al^+(16)\beta(n-32)=417781500\, \beta(n-40)>\\
77 \al^+(15)\al^+(14)  \beta(n-40)=385 945560\, \beta(n-40)>
539\al^+(16)\al^+(16)\beta(n-48)= 364193676\, \beta(n-40).$

For $n \equiv 2 \mod 4$ we have
$\al^-(16)\al^+(14)\beta(n-30)=    933147000\, \beta(n-42)>\\
77\al^+(15) \al^+(16)\beta(n-42)=863963100\, \beta(n-42)>
539\al^+(16)\al^+(14)\beta(n-46)=813454488\, \beta(n-42)$. This completes the proof of the main  statement.

Finally, by Lemma \ref{2u2}, the bound in Theorem $\ref{t21}$ is attained, proving the additional assertion. $\Box$

\medskip {\it Proof of Theorems} \ref{t21} and \ref{un5}.
Theorem \ref{t21} follows from Theorems \ref{sp88}, \ref{od1}, \ref{ev+} and \ref{ev-} 
for $q$ odd. For $q$ even use Theorems \ref{sp2}, \ref{ore2}; note that the result for the group $G^*=Spin_{2n+1}(q)$, $q$ even (not considered in  Theorems \ref{ore2}) are identical to those for $G^*=Sp_{2n}(q)$ due to the comments after Lemma \ref{f11}. Theorem \ref{un5}
follows from Theorem \ref{t21} by elementary straightforward computations.  $\Box$

\medskip

{\bf Acknowledgement.} We are very grateful to Gunter Malle for his comments on
the original manuscript which were helpful in correcting inaccuracies and improving
upon the presentation.

\bigskip

\noindent Christine Bessenrodt:
\\
Faculty of Mathematics and Physics,
\\
Leibniz University Hannover,
\\
Welfengarten 1, D-30167, Hannover,
\\
Germany

\medskip

\noindent Alexandre Zalesski:
\\
Department of Physics, Mathematics and Informatics,
\\
National Academy  of Sciences of Belarus,
\\
66 Nezalejnasti prospekt, 220072, Minsk,
\\
Belarus

\section{Appendix: The numerical data}


\begin{center}

\centerline{Table 2: $\al(n),\al^+(n)$ and $\al^-(n)$ for $1\le n\leq 43$ }\label{t1}

\bigskip
\begin{tabular}{|l|l|l|l|l|}\hline
$n$&$\beta(n)$&$\al(n)$&$\al^+(n)$&$\al^-(n)$\\
\hline
\hline
$1$&$1$ &$2$&1&1\\
$2$&$2$&$6$&4&2\\
$3$&$3$&$12$&5&5\\
$4$&$5$&$25$&14&10\\
$5$&$7$&$46$&20&20\\
$6$&$11$&$86$&42&36\\
$7$&$15$&$148$&65&65\\
$8$&$25$&$255$&120&110\\
$9$&$35$&$420$&186&186\\
$10$&$55$&$686$&316&302\\
$11$&$77$&$1088$&486&486\\
$12$&$125$&$1712$&784&762\\
$13$&$175$&$2634$&1185&1185\\
$14$&$275$&$4020$&1836&1806\\
$15$&$385$&$6036$&2730&2730\\
$16$&$625$&$8988$&4110&4066\\
$17$&$875$&$13214$&6007&6007\\
$18$&$1375$&$19282$&8830&8770\\
$19$&$1925$ &27840&12711&12711\\
$20$&3125&$39923$ &18326&18242 \\
$21$&4375&$56750$ &26007&26007\\
$22$&6875&$86160$ &36884&36772\\
$23$&9625&$112384$ &51675&51675\\
$24$&15625&$156660$ &72260&72106\\
$25$&21875&$216958$ &100058&100058\\
$26$&34375&$298894$ &138186&137984\\
$27$&48125&$409420$ &189322&189322\\
$28$&78125&$558119$ &258610&258340\\
$29$&109375&$756950$ &350877&350877\\
$30$&171875&$1022090$ &474580&474228\\
$31$&240625&$1373760$ &638203&638203\\
$32$&390625&$1838932$ &855536&855074\\
$33$&546875&$2451366$ &1141125&1141125\\
$34$&859375&$3255480$ &1517336&1516742\\
$35$&1203125&$4306920$ &2008633&2008633\\
$36$&1953125&$5678104$ &2651020&2650250\\
$37$&2734375&$7459634$ &3484969&3484969\\
$38$&4296875&$9768386$ &4568010&4567030\\
$39$&6015625&$12750360$ &5966183&5966183\\
$40$&9765625&$16592332$ &7770754&7769500\\
$41$&13671875&$21527228$ &10088066 &10088066 \\
$42$&21484375&$27850932$ &13061880&13060296\\
$43$&30078125&$35931532$ &16861595&16861595\\
\hline
\end{tabular}
\end{center}

\newpage

\begin{center} {Table 3: Maxima of $\al(a)\al(b)\beta(n-a-b)$ and $\al(a)\al ^+(b)\beta(n-a-b)$ for $1\le n \le 33$}
 \label{t2}

\bigskip
\begin{tabular}{|l|l|l|}\hline
$n$&$\al(a)\al(b)\beta(n-a-b)$&$\al(a)\al^+(b)\beta(n-a-b)$ \\
\hline
\hline
$1$ & $\al(1)\al(0)=2$ &$\al(1)\al^+(0)=2$\\
\hline
$2$ & $\al(2)\al(0)=6$ &$\al(2)\al^+(0)=6$\\
\hline
$3$ & $\al(2)\al(1)=\al(3)\al(0)=12$ &$\al(3)\al^+(0)=12$\\
\hline
$4$ & $\al(2)\al(2)=36$ &$\al(4)\al^+(0)=25$\\
\hline
$5$ & $\al(3)\al(2)=72$ &$\al(3)\al^+(2)=48 $\\
\hline
$6$ & $\al(4)\al(2)=150$ &$\al(4)\al^+(2)=100$\\
\hline
$7$ & $\al(4)\al(3)=300$ &$\al(5)\al^+(2)=184$\\
\hline
$8$ & $\al(4)\al(4)=625$ &$\al(4)\al^+(4)=350$\\
\hline
$9$ & $\al(5)\al(4)=1150$ &$\al(5)\al^+(4)=644$\\
\hline
$10$ & $\al(6)\al(4)=2150$ &$\al(6)\al^+(4)=1204$\\
\hline
$11$ & $\al(6)\al(5)=3956$ &$\al(7)\al^+(4)=2072$\\
\hline
$12$ & $\al(6)\al(6)=7396$ &$\al(6)\al^+(6)=3612$\\
\hline
$13$ & $\al(7)\al(6)=12728$ &$\al(7)\al^+(6)=6216$\\
\hline
$14$ & $\al(8)\al(6)=21930$ &$\al(8)\al^+(6)=10710$\\
\hline
$15$ & $\al(8)\al(7)=37740$ &$\al(7)\al^+(8)=17760$\\
\hline
$16$ & $\al(8)\al(8)= 65025$ &$\al(8)\al^+(8)=30600$\\
\hline
$17$ & $\al(9)\al(8)= 107100$ &$\al(9)\al^+(8)=50400$\\
\hline
$18$ & $\al(9)\al(9)=176400 $ &$\al(10)\al^+(8)=82320$\\
\hline
$19$ & $\al(10)\al(9)=288120$ &$\al(9)\al^+(10)=132720$\\
\hline
$20$ & $\al(10)\al(10)=470596 $ &$\al(10)\al^+(10)=216776$\\
\hline
$21$ & $\al(11)\al(10)=746368 $ &$\al(11)\al^+(10)=343808$\\
\hline
$22$ & $\al(11)\al(11)=1183744 $ &$\al(12)\al^+(10)=540992$\\
\hline
$23$ & $\al(12)\al(11)= 1862656$ &$\al(11)\al^+(12)=852992$\\
\hline
$24$ & $\al(12)\al(12)=2930944$ &$\al(12)\al^+(12)=1342208$\\
\hline
$25$ & $\al(13)\al(12)=4509408$ &$\al(13)\al^+(12)=2065056$\\
\hline
$26$ & $\al(13)\al(13)=6937956 $ &$\al(14)\al^+(12)=3151680$\\
\hline
$27$ & $\al(14)\al(13)=10588680$ &$\al(13)\al^+(14)=4836024$\\
\hline
$28$ & $\al(14)\al(14)=16160400$ &$\al(14)\al^+(14)=7380720$\\
\hline
$29$ & $\al(15)\al(14)=24264720  $&$\al(15)\al^+(14)=11082096$\\
\hline
$30$ & $\al(15)\al(15)=36433296$ &$\al(14)\al^+(16)=16522200$\\
\hline
$31$ & $\al(16)\al(15)=54251568$ &$\al(15)\al^+(16)=24807960$\\
\hline
$32$ & $\al(14)\al(14)\beta(4)=80802000$ &$\al(16)\al^+(16)=36940680$\\
\hline
$33$ & $\al(15)\al(14)\beta(4)=121323600$ &$\al(15)\al^+(14)\be(4)=55410480$\\
\hline
\end{tabular}
\end{center}

\bigskip

{\it Remark.} For $n\le 31$, the table also implies the maxima of the function $\al(a)\al(b)$, for $a,b$ such that $a+b=n$.
For $n=32$ the maximum is attained at $(a,b)=(16,16)$ with value 80784144;
for $n=33$ the maximum is attained at $(a,b)=(16,17)$ with value 118767432.

\newpage

\begin{center} {Table 4: Maxima of $\al^+(a)\al^+(b)\beta(n-a-b)$ and $\al^-(a)\al ^+(b)\beta(n-a-b)$ for  $1\le n \le 33$}\label{t3}

\bigskip
\begin{tabular}{|l|l|l|}\hline
$n$&$\al^+(a)\al^+(b)\beta(n-a-b)$&$\al^-(a)\al^+(b)\beta(n-a-b)$\\
\hline
\hline
$1$ & $\al^+(1)\al^+(0)=1$ &$\al^-(1)\al^+(0)=1$\\
\hline
$2$ & $\al^+(2)\al^+(0)=4$ &$\al^-(0)\al^+(2)=4$\\
\hline
$3$ & $\al^+(3)\al^+(0)=5$ &$\al^-(3)\al^+(0)=5$\\
\hline
$4$ & $\al^+(2)\al^+(2)=16$ &$\al^-(0)\al^+(4)=14$\\
\hline
$5$ & $\al^+(3)\al^+(2)=\al^+(5)\al^+(0)=20$
&$\al^-(3)\al^+(2)=\al^-(5)\al^+(0)=20$\\
\hline
$6$ & $\al^+(4)\al^+(2)=56$ &$\al^-(0)\al^+(6)=42$\\
\hline
$7$ & $\al^+(5)\al^+(2)=80$ &$\al^-(5)\al^+(2)=80$\\
\hline
$8$ & $\al^+(4)\al^+(4)=196$ &$\al^-(6)\al^+(2)=144$\\
\hline
$9$ & $\al^+(5)\al^+(4)=280$ &$\al^-(5)\al^+(4)=280$\\
\hline
$10$ & $\al^+(6)\al^+(4)=588$ &$\al^-(6)\al^+(4)=504$\\
\hline
$11$ & $\al^+(7)\al^+(4)=910$ &$ \al^-(7)\al^+(4)=910$\\
\hline
$12$ & $\al^+(6)\al^+(6)=1764$ &$\al^-(8)\al^+(4)=1540$\\
\hline
$13$ & $\al^+(7)\al^+(6)=2730$ &$\al^-(7)\al^+(6)=2730$\\
\hline
$14$ & $\al^+(8)\al^+(6)=5040$ &$\al^-(8)\al^+(6)=4620$\\
\hline
$15$ & $\al^+(9)\al^+(6)=7812$ &$\al^-(9)\al^+(6)=7812$\\
\hline
$16$ & $\al^+(8)\al^+(8)=14400$ &$\al^-(8)\al^+(8)=13200$\\
\hline
$17$ & $\al^+(9)\al^+(8)=22320$ &$\al^-(9)\al^+(8)=22320$\\
\hline
$18$ & $\al^+(10)\al^+(8)=37920$ &$\al^-(10)\al^+(8)=36240$\\
\hline
$19$ & $\al^+(9)\al^+(10)=58776$ &$\al^-(9)\al^+(10)=58776$\\
\hline
$20$ & $\al^+(10)\al^+(10)=99856$ &$\al^-(10)\al^+(10)=95432$\\
\hline
$21$ & $\al^+(11)\al^+(10)=153576$ &$\al^-(11)\al^+(10)=153576$\\
\hline
$22$ & $\al^+(12)\al^+(10)=247744$ &$\al^-(12)\al^+(10)=240792$\\
\hline
$23$ & $\al^+(11)\al^+(12)=381024$ &$\al^-(11)\al^+(12)=381024$\\
\hline
$24$ & $\al^+(12)\al^+(12)=614656$ &$\al^-(12)\al^+(12)=397408$\\
\hline
$25$ & $\al^+(13)\al^+(12)=929040$ &$\al^-(13)\al^+(12)=929040$\\
\hline
$26$ & $\al^+(14)\al^+(12)=1439424$ &$\al^-(14)\al^+(12)=1415904$\\
\hline
$27$ & $\al^+(13)\al^+(14)=2175660$ &$\al^-(13)\al^+(14)=2175660$\\
\hline
$28$ & $\al^+(14)\al^+(14)=3370896$ &$\al^-(14)\al^+(14)=3315816$\\
\hline
$29$ & $\al^+(15)\al^+(14)= 5012280$&$\al^-(15)\al^+(14)=5012280$\\
\hline
$30$ & $\al^+(16)\al^+(14)=7545960  $ &$\al^-(16)\al^+(14)=7465176$\\
\hline
$31$ & $\al^+(15)\al^+(16)=11220300$ &$\al^-(15)\al^+(16)=11220300$\\
\hline
$32$ & $\al^+(16)\al^+(16)=16892100$ &$\al^-(16)\al^+(16)=16711260$\\
\hline
$33$ & $\al^+(15)\al^+(14)\beta(4)=25061400$ &$\al^-(15)\al^+(14)\beta(4)=25061400$\\
\hline
\end{tabular}
\end{center}

\end{document}